\documentclass[11pt,reqno]{amsart}
\usepackage[T1]{fontenc}
\usepackage{lmodern}
\usepackage{amsmath,amssymb,amsfonts,amsthm,mathtools}
\usepackage{enumitem}
\usepackage{microtype}
\usepackage{booktabs}
\usepackage{array}
\usepackage[colorlinks=true,linkcolor=blue,citecolor=blue,urlcolor=blue]{hyperref}
\usepackage[margin=0.84in]{geometry}
\usepackage[capitalise,nameinlink,noabbrev]{cleveref}
\allowdisplaybreaks
\newtheorem{theorem}{Theorem}[section]
\newtheorem{proposition}[theorem]{Proposition}
\newtheorem{lemma}[theorem]{Lemma}
\newtheorem{corollary}[theorem]{Corollary}
\theoremstyle{definition}
\newtheorem{definition}[theorem]{Definition}
\newtheorem{remark}[theorem]{Remark}

\numberwithin{equation}{section}

\newcommand{\R}{\mathbb R}
\newcommand{\F}{\mathcal F}
\newcommand{\Lop}{\mathcal L}
\newcommand{\Sop}{\mathcal S}
\newcommand{\CT}{\mathcal C}
\newcommand{\e}{\mathrm e}
\newcommand{\polhk}[1]{\k{#1}}

\title[Doubly critical semilinear heat equations]{Nonuniqueness and Failure of Positivity Preservation for Doubly Critical
Semilinear Heat Equations}
\author{Berikbol T. Torebek}
\address{Institute of Mathematics and Mathematical Modeling,
28 Shevchenko str., 050010 Almaty, Kazakhstan; and
SDU University, Kaskelen, Kazakhstan}
\email{torebek@math.kz}
\date{}

\begin{document}

\begin{abstract}
At the Fujita exponent $p=1+2/N$ we construct a strictly positive mild
solution of $u_t-\Delta u=|u|^{p-1}u$ in $\mathbb R^N$ whose initial trace in
$L^1$ is zero. In similarity variables the solution approaches a Gaussian
profile, and
\[
 \|U(t)\|_1\sim C_N(\log(1/t))^{-N/2},\qquad
 \|U(t)\|_\infty\sim c_Nt^{-N/2}(\log(1/t))^{-N/2},
\]
with explicit constants $C_N$ and $c_N$. In particular, uniqueness fails at
the borderline logarithmic rate, and the solution lies outside every stronger
logarithmic uniqueness class; this answers two questions raised by Miyamoto
{\it [J. Evol. Equ. 21 (2021), 151--166, Remark~1.6]}. We also prove that, up to
time translation, this is the only ancient solution of fixed sign of the
rescaled equation that decays in a Gaussian weighted supremum norm.

On a bounded smooth domain the whole-space solution can be localized at
finitely many interior points, with independently chosen signs, on top of an
arbitrary bounded initial datum. The resulting solutions are classical for
$t>0$ and continuous at $t=0$ in $L^1$. Consequently, uniqueness fails for
every bounded datum, and positivity is lost for every bounded nonnegative
datum, $u_0\equiv1$ included. These results provide affirmative answers to Open Problems~3 and~4 posed by Brezis and Cazenave in {\it [J. Anal. Math. 68 (1996), 277--304]}.
\end{abstract}



\maketitle

\section{Introduction}\label{sec:intro}

Let $N\ge1$, $0<\alpha<1$, and let $\Omega\subset\R^N$ be a bounded domain
with $C^{2,\alpha}$ boundary. We consider the Cauchy--Dirichlet problem
\begin{equation}\label{eq:main}
\begin{cases}
 u_t-\Delta u=|u|^{p-1}u,&(t,x)\in(0,T)\times\Omega,\\
 u=0,&(t,x)\in(0,T)\times\partial\Omega,\\
 u(0)=u_0,&x\in\Omega,
\end{cases}
\end{equation}
as well as the Cauchy problem in $\R^N$. For a power nonlinearity the
relevant Lebesgue exponent is $q_c=N(p-1)/2$, for which the $L^{q_c}$ norm is
invariant under the scaling of the equation. When $q>q_c$, or $q=q_c>1$, the
semigroup approach of Weissler \cite{Weissler1980} and the analysis of Brezis
and Cazenave \cite{BrezisCazenave1996} give local solutions with the usual
smoothing properties, while for $q<q_c$ local solvability may fail.
Unconditional uniqueness in $C([0,T];L^q)$ is a more delicate matter, and at
the critical exponent it has to be treated separately. At the Fujita exponent
\begin{equation*}
 p=p_F:=1+\frac2N
\end{equation*}
one has $q_c=1$, so that $L^1$ itself is the critical space; Brezis and
Cazenave called this the \emph{doubly critical} case. Here the $L^1$ theory
behaves differently, and several of the questions listed in
\cite[Section~7.5]{BrezisCazenave1996} have remained open.

Two of these questions concern the class
\begin{equation*}
 \CT_T(\Omega):=C([0,T];L^1(\Omega))\cap L^\infty_{\rm loc}((0,T);L^\infty(\Omega)).
\end{equation*}
Brezis and Cazenave asked whether uniqueness in $\CT_T(\Omega)$ can fail for
some $u_0\in L^1(\Omega)$ (Open Problem~3), and whether a solution starting from
$u_0\ge0$ can become negative (Open Problem~4). A third question from the same
list, namely whether there is an $L^1$ datum for which \eqref{eq:main} has no
local solution at all, was answered by Celik and Zhou \cite{CelikZhou2003}; see
also \cite{LRSV2016} and \cite[Theorem~1.3]{Miyamoto2021}. In the formulation of
Brezis and Cazenave the datum $u_0=0$ is allowed in both problems, and we shall
use this. In fact we show that both phenomena occur for every bounded datum.
For the positivity question the datum may even be strictly positive: it need
not vanish near the point where the negative part concentrates.

Already on $\R^N$ the endpoint case is subtle. Miyamoto \cite{Miyamoto2021}
proved local existence for data in the logarithmic Orlicz class $X_{N/2}$, where
$X_\ell=\{\phi:\int_{\R^N}|\phi|[\log(\e+|\phi|)]^\ell<\infty\}$, and showed
that the exponent $N/2$ is sharp in this scale: for every $0\le\ell<N/2$ there
are nonnegative data in $X_\ell$ without nonnegative local solutions (the proof
relies on the criterion of Baras and Pierre \cite{BarasPierre1985}). As for
uniqueness, he proved that an integral solution is unique in
\begin{equation*}
 \Bigl\{u:\ \sup_{0<t<T}t^{N/2}(-\log t)^\ell\|u(t)\|_\infty<\infty\Bigr\}
\end{equation*}
whenever $\ell>N/2$, and he pointed out \cite[Remark~1.6]{Miyamoto2021} that
two questions remained open: can a solution exist outside all of these classes,
and does uniqueness persist at the endpoint $\ell=N/2$? We also mention the work
of Laister and Sier\.z\polhk{e}ga \cite{LaisterSierzega2020}, who obtained $L^1$
well-posedness in $\CT_T$ for general nonlinearities under an integrability
condition which, for the pure power, holds exactly when $p<p_F$. The Fujita
exponent lies just outside their criterion, in agreement with the
nonuniqueness established below.

None of the known nonuniqueness mechanisms reaches the endpoint $p=p_F$. Above
the Fujita exponent, positive self-similar solutions with zero initial trace go
back to Haraux and Weissler \cite{HarauxWeissler1982}, and their profiles were
studied in
\cite{Weissler1986a,Weissler1986b,EscobedoKavian1987,Yanagida1996,DohmenHirose1998};
see also Baras \cite{Baras1983}. More recently, Glogi\'c, Hofmanov\'a, Lange and
Luongo \cite{GlogicHofmanovaLangeLuongo2025} proved nonuniqueness in
$L^q(\R^N)$, $1\le q<q_c$, for $N\ge3$ and $p_F<p<p_{JL}$. All these
constructions need a nonempty range of Lebesgue exponents $q<q_c$, and this
range is empty when $p=p_F$. Other mechanisms start from singular initial data
or singular stationary solutions
\cite{Naito2004,Naito2006,SoupletWeissler2003,CazenaveDicksteinNaumkinWeissler2020,CazenaveDicksteinNaumkinWeissler2019,FujishimaIokuRufTerraneo2026,HisaMiyamoto2026},
or use connecting orbits and ancient solutions of a rescaled equation
\cite{FilaYanagida2011,FilaMatanoYanagida2017}; we compare these results with
ours in \cref{sec:literature}. The new ingredient here is a logarithmic
modulation of the self-similar amplitude, which takes over the role of the
nonzero self-similar profile once the subcritical Lebesgue range disappears.

\paragraph{Notion of solution.}
We write $\F(s):=|s|^{p-1}s$, and we denote by $e^{t\Delta}$ the heat semigroup
on $\R^N$ and by $S_\Omega(t)$ the Dirichlet heat semigroup on $\Omega$.

\begin{definition}\label{def:sol}
Let $u_0\in L^1(\Omega)$ and $T>0$. A \emph{solution} of \eqref{eq:main} on
$[0,T]$ is a function $u\in\CT_T(\Omega)$ that is a classical solution of the
equation on $(0,T)\times\Omega$, is continuous on
$(0,T)\times\overline\Omega$ with $u=0$ on $(0,T)\times\partial\Omega$, and
satisfies $u(0)=u_0$.
\end{definition}

This is the class in which Open Problems~3 and~4 are formulated; Brezis and
Cazenave consider solutions that are classical for $t>0$ and take the initial
value in $L^q$. All solutions constructed below satisfy, in addition, the
Duhamel formula
\begin{equation}\label{eq:mildOmega}
 u(t)=S_\Omega(t)u_0+\int_0^tS_\Omega(t-s)\F(u(s))\,ds,\qquad
 \F(u)\in L^1((0,T);L^1(\Omega)).
\end{equation}
On $\R^N$ we use the notion of Miyamoto \cite[Definition~1.2]{Miyamoto2021}: a
\emph{mild solution} with datum $\phi\in L^1(\R^N)$ is a function
$u\in C([0,T);L^1(\R^N))$ with $\|u(t)\|_\infty<\infty$ for $0<t<T$ such that
\[
u(t)=e^{t\Delta}\phi+\int_0^te^{(t-s)\Delta}\F(u(s))\,ds.
\]

\subsection{Main results}

Our results concern the construction on the whole space, its rigidity in
Gaussian weighted spaces, and its localization on bounded domains. Throughout we
use the constants
\begin{equation*}
 c_N:=\Bigl(\frac N2\Bigr)^{N/2}\Bigl(1+\frac2N\Bigr)^{N^2/4},
 \qquad
 C_N:=(4\pi)^{N/2}c_N
 =(2\pi N)^{N/2}\Bigl(1+\frac2N\Bigr)^{N^2/4},
\end{equation*}
and we write $L(t):=\log(1/t)$ and $T_*:=\e^{-2}$.

\begin{theorem}[Positive whole-space zero-trace branch]\label{thm:whole}
Let $N\ge1$ and $p=1+2/N$.  There is a classical solution $U$ of
\begin{equation}\label{eq:whole}
 U_t-\Delta U=|U|^{p-1}U\qquad\text{in }(0,1]\times\R^N
\end{equation}
such that
\[
 U\in C([0,1];L^1(\R^N))\cap L^\infty_{\rm loc}((0,1];L^\infty(\R^N)),
 \qquad U(0)=0,
\]
and $U$ is a mild solution with datum $0$ on $[0,1]$.  There are constants
$c,C>0$ such that, for $0<t\le T_*$,
\begin{equation}\label{eq:wholebounds}
\begin{aligned}
 cL(t)^{-N/2}\le\|U(t)\|_1&\le CL(t)^{-N/2},\\
 \|U(t)\|_\infty&\le Ct^{-N/2}L(t)^{-N/2},\\
 |U(t,x)|&\le Ct^{-N/2}L(t)^{-N/2}
 \exp\!\Bigl(-c\frac{|x|^2}{t}\Bigr).
\end{aligned}
\end{equation}
As $t\downarrow0$,
\begin{equation}\label{eq:exactnorms}
 \|U(t)\|_1\sim C_NL(t)^{-N/2},\qquad
 \|U(t)\|_\infty\sim c_Nt^{-N/2}L(t)^{-N/2},\qquad
 t^{N/2}L(t)^{N/2}U(t,0)\to c_N .
\end{equation}
Moreover $U(t,x)>0$ for $0<t\le1$ and $x\in\R^N$. Consequently
\begin{equation}\label{eq:Miy}
 \sup_{0<t<T_*}t^{N/2}L(t)^{N/2}\|U(t)\|_\infty<\infty,
 \qquad
 t^{N/2}L(t)^\ell\|U(t)\|_\infty\longrightarrow\infty
 \quad(t\downarrow0,\ \ell>N/2).
\end{equation}
\end{theorem}

\begin{remark}[Miyamoto's questions]\label{rem:Miyamoto}
The first estimate in \eqref{eq:Miy} says that $U$ belongs to the endpoint
class $\ell=N/2$, and the second that it belongs to no class with $\ell>N/2$.
Thus there is a nontrivial solution in the sharp logarithmic class, but in no
smaller one. Since $u\equiv0$ is another nonnegative solution with the same
datum, uniqueness fails at the endpoint even among nonnegative solutions. This
answers both questions in \cite[Remark~1.6]{Miyamoto2021}. We also note that
$U$ attains, up to constants, the upper bound of Hisa and Ishige
\cite{HisaIshige2018} for the mass of nonnegative solutions in parabolic
balls. Thus the rate $L(t)^{-N/2}$ in \eqref{eq:wholebounds} is the largest
possible for this mass among nonnegative solutions; see
\cref{sec:absorption}.
\end{remark}

\begin{corollary}[Scaling family]\label{cor:scaling}
Let $U$ be the solution of \cref{thm:whole}.  For every $\lambda\in(0,1]$ the
rescaling
\[
 U_\lambda(t,x):=\lambda^N U(\lambda^2t,\lambda x)
\]
is a positive mild solution of \eqref{eq:whole} on $(0,1]$ with zero initial
trace, and the family $\{U_\lambda:\lambda\in(0,1]\}$ is pairwise distinct.
\end{corollary}

The spaces $X_\theta$, the Gaussian $\rho$ and the operator $\Lop$ are
introduced in \cref{sec:sim}. The next theorem is a rigidity statement for
ancient solutions in these weighted spaces; we stress that it does not classify
all solutions with zero initial trace in $L^1$.

\begin{theorem}[Fixed-sign rigidity in a Gaussian space]\label{thm:rigidity}
Let $0<\theta<1$ and let
$v\in C((-\infty,0];X_\theta)$ be a classical solution of
\[
 v_\tau=\Lop v+|v|^{p-1}v\qquad\text{in }(-\infty,0)\times\R^N
\]
such that $\|v(\tau)\|_\theta\to0$ as $\tau\to-\infty$.  Assume that
$v\not\equiv0$ and that $v$ has one fixed sign.  Then, after multiplication by
$-1$ if necessary, the following hold.
\begin{enumerate}[label=\textup{(\roman*)}]
\item $v>0$ on $(-\infty,0]\times\R^N$, and
\[
 v(\tau)=a(\tau)\rho+z(\tau),\qquad a(\tau)=\int_{\R^N}v(\tau,y)\,dy,\qquad
 \|z(\tau)\|_\theta=O(a(\tau)^p),
\]
with
\begin{equation*}
 a(\tau)=C_N(-\tau)^{-N/2}
 \left[1+O\!\left(\frac{\log(-\tau)}{-\tau}\right)\right]
 \qquad(\tau\to-\infty).
\end{equation*}
Consequently
\begin{equation*}
 (-\tau)^{N/2}v(\tau,\cdot)\longrightarrow C_N\rho
 \qquad\text{in }X_\theta.
\end{equation*}
\item Let $\theta_0,\theta_1$ be chosen as in \eqref{eq:weightschoice}. There
is $\bar\tau<0$ such that, for every $\tau_*\le\bar\tau$, the
translate $v(\cdot+\tau_*)$ restricted to $(-\infty,0]$ is the
Lyapunov--Perron orbit $v_\varepsilon$ of \cref{prop:LP} with terminal mass
$\varepsilon=a(\tau_*)$.
\item For every admissible $\bar\varepsilon$ in \cref{prop:LP} there is
$s\in\R$ such that $v(\tau)=v_{\bar\varepsilon}(\tau+s)$ for all
$\tau\le\min\{0,-s\}$.
\end{enumerate}
\end{theorem}

\begin{remark}[Uniqueness up to scaling]\label{rem:rigidscaling}
A time translation $\tau\mapsto\tau+s$ in similarity variables corresponds to
the scaling $U\mapsto\lambda^NU(\lambda^2t,\lambda x)$ with $\lambda=\e^{s/2}$.
Part (iii) of \cref{thm:rigidity} therefore says that, within this Gaussian class, the solution of
\cref{thm:whole} is unique up to the rescalings of \cref{cor:scaling} (now with
arbitrary $\lambda>0$, on the common interval of definition). The choice of the
weight is irrelevant: since $X_\theta\subset X_{\theta'}$ for $\theta'<\theta$,
part (iii) applied with the smaller of two weights shows that orbits constructed
with different admissible weights, and in particular the similarity profile of
$U$, are time translates of each other.
\end{remark}

\begin{remark}[The original variables]\label{rem:rigidphysical}
In the original variables \cref{thm:rigidity} reads as follows. Suppose that $u$ is a
classical solution of \eqref{eq:whole} on $(0,1]\times\R^N$, and that for some
$x_0\in\R^N$ the function $v(\tau,y):=\e^{N\tau/2}u(\e^\tau,x_0+\e^{\tau/2}y)$
satisfies the hypotheses of \cref{thm:rigidity}. Then
\begin{equation*}
 t^{N/2}L(t)^{N/2}
 u(t,x_0+\sqrt t\,\cdot)\longrightarrow \pm C_N\rho
 \qquad\text{in }X_\theta\quad(t\downarrow0),
\end{equation*}
and $u(t,x_0+\cdot)=\pm U_\lambda(t,\cdot)$ for some $\lambda>0$ and all
small $t>0$. In other words,
\[
 u(t,x_0+\sqrt t\,y)
 =\pm C_Nt^{-N/2}L(t)^{-N/2}\rho(y)
 +o_{X_\theta}\!\bigl(t^{-N/2}L(t)^{-N/2}\bigr).
\]
Here the error is measured in $X_\theta$, whose weight is weaker than the
Gaussian decay of $\rho$, so this does not give a uniform relative error in the
far field.
\end{remark}

\begin{theorem}[Arbitrary bounded data and prescribed finite signed cores]
\label{thm:multicore}
Let $\Omega\subset\R^N$ be a bounded domain with $C^{2,\alpha}$ boundary,
$p=1+2/N$, and let $u_0\in L^\infty(\Omega)$ be arbitrary.  Let
$x_1,\dots,x_m\in\Omega$ be distinct and let
$\boldsymbol\sigma=(\sigma_1,\dots,\sigma_m)\in\{-1,+1\}^m$.  Then there are
$T>0$ and two distinct solutions $V$ and $u_{\mathbf x,\boldsymbol\sigma}$ of
\eqref{eq:main} on $[0,T]$, in the sense of \cref{def:sol} and satisfying
\eqref{eq:mildOmega}, with the following properties.
\begin{enumerate}[label=\textup{(\roman*)}]
\item $V\in L^\infty((0,T)\times\Omega)$ is the ordinary bounded solution.  If
$u_0\ge0$, then $V\ge0$.
\item There is $C>0$ such that
\begin{align*}
 \|u_{\mathbf x,\boldsymbol\sigma}(t)-V(t)\|_{L^\infty(\Omega)}
 &\le Ct^{-N/2}L(t)^{-N/2},\\
 \|u_{\mathbf x,\boldsymbol\sigma}(t)-V(t)\|_{L^1(\Omega)}
 &\le CL(t)^{-N/2}.
\end{align*}
\item At every prescribed core,
\begin{equation}\label{eq:mastercenter}
 u_{\mathbf x,\boldsymbol\sigma}(t,x_j)-V(t,x_j)
 =\sigma_j c_Nt^{-N/2}L(t)^{-N/2}(1+o(1))
 \qquad(t\downarrow0),
\end{equation}
for $j=1,\dots,m$.
\item If $u_0\ge0$ and $\sigma_j=-1$ for some $j$, then
$u_{\mathbf x,\boldsymbol\sigma}(t,x_j)<0$ for all sufficiently small $t>0$.
\end{enumerate}
In particular, uniqueness in $\CT_T(\Omega)$ fails for every bounded datum, and
positivity preservation fails for every bounded nonnegative datum.  Moreover,
for every compact set $\mathcal K\Subset\Omega$ there is a common time
$T_{\mathcal K}>0$, depending also on $u_0$, such that the single-core
solutions $u_{x,\sigma}$, $x\in\mathcal K$, $\sigma\in\{-1,+1\}$, are all
defined on $[0,T_{\mathcal K}]$ and are pairwise distinct.  Hence every bounded
datum admits an uncountable family of solutions in the critical class.
\end{theorem}

\begin{remark}[Open Problems~3 and~4]\label{rem:openproblems}
Taking a single core in \cref{thm:multicore}, we obtain two distinct solutions
in $\CT_T(\Omega)$ with the same bounded datum, which answers Open Problem~3 of
Brezis and Cazenave \cite[Section~7.5]{BrezisCazenave1996}. If $u_0\ge0$ and the
core is negative, the solution takes negative values at the center of the core
for all sufficiently small $t>0$, which answers Open Problem~4 of
\cite[Section~7.5]{BrezisCazenave1996}; this happens even for strictly positive
data such as $u_0\equiv1$. The loss of positivity is confined to the initial
time: for each $t_0>0$ the solution is classical on $[t_0,T]\times\Omega$, so
the maximum principle applies there, while the $L^1$ norm of the negative core
tends to zero as $t\downarrow0$, so the core is not seen by the convergence
$u(t)\to u_0$ in $L^1$; see \cref{sec:posmass}.
\end{remark}

\subsection{Ideas of the proofs}

The construction is most transparent in similarity variables. There the
Gaussian spans the kernel of the Ornstein--Uhlenbeck operator, and the nonlinear
dynamics along this neutral direction is slow. To leading order the amplitude
$a$ of the Gaussian mode satisfies $a'=\kappa_Na^p$ with $\kappa_N=\int\rho^p$,
and the solutions of this equation with $a(-\infty)=0$ behave like
\[
 a(\tau)\sim C_N(-\tau)^{-1/(p-1)}=C_N(-\tau)^{-N/2}.
\]
Since $-\tau=\log(1/t)$, the factor $L(t)^{-N/2}$ is dictated by the neutral
mode; it is not a by-product of the method. We obtain the corresponding ancient
orbit by a Lyapunov--Perron argument. The relative error
$O(\log|\tau|/|\tau|)$ in \eqref{eq:arelerr} comes from integrating
$B-\kappa_N=O(1/|\tau|)$; we have not tried to decide whether it is sharp. The
same amplitude law, with the opposite sign and forward in time, governs the
absorption equation at the Fujita exponent (see \cref{sec:absorption}).
Positivity is obtained from approximations on finite time intervals and a
comparison argument in the mild formulation. For the rigidity result, a
spectral gap in a fixed weighted space forces every small ancient solution of
fixed sign onto the same center orbit.

On a bounded domain we cut off translated copies of the whole-space solution in
disjoint balls. The cutoff error is exponentially small as $t\downarrow0$, and
the critical term is the interaction between a bounded background and a
core:
\[
 |W(t,x)|^{p-1}\lesssim \frac1{tL(t)}
 e^{-c|x-x_j|^2/t}.
\]
The factor $1/(tL(t))$ is not integrable in time, but the Gaussian factor
compensates for this:
\[
 \left\|\int_0^tS_\Omega(t-s)
 \frac{e^{-c|\cdot-x_j|^2/s}}{sL(s)}\,ds\right\|_\infty
 \lesssim L(t)^{-1}.
\]
This is enough to close a contraction argument in the norm
$\sup_{0<t\le T}L(t)\|w(t)\|_\infty$. The same Duhamel equation gives a better
$L^1$ bound for the correction, so that the singular part keeps the mass scale
$L(t)^{-N/2}$ of the whole-space solution.

\paragraph{Organization of the paper.}
\Cref{sec:wholespace} deals with the whole space. After collecting estimates for
the Ornstein--Uhlenbeck semigroup in Gaussian weighted spaces, we construct the
ancient orbit, prove that it is positive, and deduce \cref{thm:whole},
\cref{cor:scaling} and \cref{thm:rigidity}. In \cref{sec:localization} we
localize the construction in bounded domains and prove \cref{thm:multicore}.
\Cref{sec:literature} contains further remarks and a comparison with earlier
work. The appendix collects a Gaussian identity, a lemma on Duhamel integrals
starting at $-\infty$, and the regularity facts used in the paper.

\section{The critical zero-trace solution on \texorpdfstring{$\mathbb{R}^N$}{R to the N}}\label{sec:wholespace}
\subsection{Similarity variables and the Ornstein--Uhlenbeck semigroup}\label{sec:sim}

Throughout,
\begin{equation*}
 p=1+\frac2N,\qquad \beta:=\frac1{p-1}=\frac N2,\qquad \F(s)=|s|^{p-1}s.
\end{equation*}
On $\R^N$ we use the similarity variables
\begin{equation*}
 \tau=\log t,\qquad y=\frac{x}{\sqrt t},\qquad
 U(t,x)=t^{-N/2}v(\tau,y).
\end{equation*}
Since $Np/2=N/2+1$, a direct computation gives
\[
U_t-\Delta U-\F(U)=t^{-N/2-1}\bigl(v_\tau-\Lop v-\F(v)\bigr),
\]
so that \eqref{eq:whole} is equivalent to
\begin{equation}\label{eq:veq}
 v_\tau=\Lop v+\F(v),\qquad
 \Lop=\Delta+\frac12y\cdot\nabla+\frac N2=\nabla\cdot\Bigl(\nabla+\frac y2\Bigr).
\end{equation}
The normalized Gaussian
\begin{equation*}
 \rho(y)=(4\pi)^{-N/2}e^{-|y|^2/4}
\end{equation*}
satisfies $\Lop\rho=0$ and $\int\rho=1$. We denote by $\Sop(s)=e^{s\Lop}$,
$s\ge0$, the Ornstein--Uhlenbeck semigroup generated by $\Lop$; it is given by
Mehler's formula
\begin{equation*}
 (\Sop(s)f)(y)=\int_{\R^N}K_\omega(y,x)f(x)\,dx,\qquad
 K_\omega(y,x)=\frac{1}{[4\pi(1-\omega^2)]^{N/2}}
 \exp\!\Bigl(-\frac{|y-\omega x|^2}{4(1-\omega^2)}\Bigr),
\end{equation*}
where $\omega=e^{-s/2}$; at $s=\infty$, i.e.\ $\omega=0$, we read
$K_0(y,x)=\rho(y)$. Clearly $\Sop$ is positive, preserves the integral and
leaves $\rho$ invariant. In $L^2(e^{|y|^2/4}dy)$ the operator $\Lop$ is
self-adjoint with eigenvalues $-k/2$, $k=0,1,2,\dots$; see
\cite{Kavian1987,EscobedoKavian1987}, and \cite{MetafunePallaraPriola2002} for
the spectrum in $L^p$ spaces with respect to the invariant measure. We need,
however, estimates in weighted supremum norms, and these we derive directly
from the kernel.

The substitution $x'=\sqrt s\,\xi$ in the heat kernel shows that, for
$h\in L^1(\R^N)$ and $0<s<t$,
\begin{equation}\label{eq:conj}
 \Bigl(e^{(t-s)\Delta}[s^{-N/2}h(\cdot/\sqrt s)]\Bigr)(x)
 =t^{-N/2}\bigl(\Sop(\log(t/s))h\bigr)(x/\sqrt t),
\end{equation}
with $\omega=\sqrt{s/t}$ in Mehler's formula. This identity allows us to pass
back and forth between the two formulations of the problem.

For $0<\theta<1$ let
\begin{equation*}
 X_\theta=\Bigl\{f\in C(\R^N):
 \|f\|_\theta:=\sup_{y\in\R^N}e^{\theta|y|^2/4}|f(y)|<\infty\Bigr\}.
\end{equation*}
Then $\|f\|_{L^r}\le(4\pi/(r\theta))^{N/(2r)}\|f\|_\theta$ for $1\le r<\infty$,
$\|f\|_\infty\le\|f\|_\theta$, and $\|f\|_{\theta'}\le\|f\|_\theta$ for
$\theta'\le\theta$.

We note that $\Sop$ is not strongly continuous on $X_\theta$. For instance,
if $f_\theta(y)=e^{-\theta|y|^2/4}$, Gaussian integration gives
\[
 \Sop(s)f_\theta(y)
 =D_s^{-N/2}e^{-\theta|y|^2/(4D_s)},\qquad
 D_s=\theta+(1-\theta)e^{-s}.
\]
For $s>0$ we have $D_s<1$, so
$e^{\theta|y|^2/4}\Sop(s)f_\theta(y)\to0$ as $|y|\to\infty$,
whereas $e^{\theta|y|^2/4}f_\theta(y)=1$. Thus
$\|\Sop(s)f_\theta-f_\theta\|_\theta\ge1$.

What we shall use instead is continuity with a loss in the weight.
If $0<\theta<\widehat\theta<1$ and $f\in X_{\widehat\theta}$, then
$s\mapsto\Sop(s)f$ is continuous from $[0,\infty)$ into $X_\theta$,
with $\Sop(0)f=f$. To see this, Gaussian integration gives the uniform
bound
\[
 e^{\theta|y|^2/4}|\Sop(s)f(y)|
 \le\widehat\theta^{-N/2}\|f\|_{\widehat\theta}
 e^{-(\widehat\theta-\theta)|y|^2/4},\qquad s\ge0.
\]
The tails are therefore uniformly small in the weaker weighted norm.
On a fixed ball, continuity follows from the Mehler kernel for $s>0$
and its approximate-identity property at $s=0$. Together with the uniform
bound from $X_{\widehat\theta}$ to $X_\theta$, this also gives continuity
when the input varies continuously in $X_{\widehat\theta}$.

Let
\begin{equation*}
 Pf=\Bigl(\int_{\R^N}f\Bigr)\rho,\qquad Q=I-P,
\end{equation*}
so that
\[
P\Sop(s)=\Sop(s)P=P\qquad \text{and}\qquad Q\Sop(s)=\Sop(s)Q.
\]
From $|\int f|\le(4\pi/\theta)^{N/2}\|f\|_\theta$ and
$\|\rho\|_\theta=(4\pi)^{-N/2}$, we obtain
\begin{equation}\label{eq:Pbound}
 \|Pf\|_\theta\le\theta^{-N/2}\|f\|_\theta,\qquad
 \|Qf\|_\theta\le(1+\theta^{-N/2})\|f\|_\theta .
\end{equation}

\begin{lemma}[Ornstein--Uhlenbeck estimates in $X_\theta$]\label{lem:OU}
Let $0<\theta\le\theta_1<1$ and $s\ge0$.
\begin{enumerate}[label=\textup{(\alph*)}]
\item For $f\in X_{\theta_1}$,
\begin{equation*}
 \|\Sop(s)f\|_\theta\le\theta_1^{-N/2}\|f\|_{\theta_1}.
\end{equation*}
\item There is $C_\theta>0$ such that, for $f\in X_\theta$,
\begin{equation}\label{eq:Sgapsame}
 \|\Sop(s)Qf\|_\theta\le C_\theta e^{-s/2}\|f\|_\theta .
\end{equation}
Consequently, with $C_{\rm gap}:=C_\theta$,
\begin{equation}\label{eq:Sgap}
 \|\Sop(s)Qf\|_\theta\le C_{\rm gap}\,e^{-s/2}\|f\|_{\theta_1}
 \qquad(f\in X_{\theta_1}).
\end{equation}
\item The rate is optimal: $\Sop(s)Q\partial_{y_1}\rho=e^{-s/2}\partial_{y_1}\rho$.
\item Let $0<\theta_0<\theta_1<1$ and
$k(s)=s^{-1/2}\mathbf 1_{(0,1]}(s)+e^{-s/2}\mathbf 1_{(1,\infty)}(s)$.
There is $C_\nabla=C_\nabla(N,\theta_0,\theta_1)$ such that
$\Sop(s)Qf\in C^1(\R^N)$ for $f\in X_{\theta_1}$, $s>0$, and
\begin{equation*}
 \|\nabla\Sop(s)Qf\|_{\theta_0}\le C_\nabla\,k(s)\,\|f\|_{\theta_1}.
\end{equation*}
\end{enumerate}
\end{lemma}

\begin{proof}
We first record a kernel bound used in (b) and (d). Let
$0<\theta'\le\theta''<1$, $\omega\in[0,e^{-1/2}]$ and $b=1-\omega^2$. The
logarithmic derivatives of the Mehler kernel are
\begin{align*}
 \partial_\omega\log K_\omega&=\frac{N\omega}{b}+\frac{(y-\omega x)\cdot x}{2b}
 -\frac{\omega|y-\omega x|^2}{2b^2},\qquad
 \nabla_y\log K_\omega=-\frac{y-\omega x}{2b},\\
 \partial_\omega\nabla_y\log K_\omega&=\frac{x}{2b}-\frac{\omega(y-\omega x)}{b^2}.
\end{align*}
Hence, for $j,k\in\{0,1\}$, $\partial_\omega^j\nabla_y^kK_\omega=K_\omega
\mathcal P_{j,k}$, where $\mathcal P_{j,k}$ is a polynomial in $(x,y)$ of
degree at most $2j+k$ with coefficients continuous on $[0,e^{-1/2}]$. For each
coordinate pair $(x_i,y_i)$,
\[
 \frac{(y_i-\omega x_i)^2}{b}+\theta''x_i^2-\theta'y_i^2
 \ge\frac{(y_i-\omega x_i)^2}{b}+\theta'x_i^2-\theta'y_i^2 .
\]
The quadratic form on the right has determinant $\theta'(1-\theta')$ and trace
$(1+\omega^2)/b\le2/(1-e^{-1})$, so its smallest eigenvalue is at least
$c_1:=\theta'(1-\theta')(1-e^{-1})/2$. Therefore
$e^{\theta'|y|^2/4}K_\omega(y,x)e^{-\theta''|x|^2/4}\le
Ce^{-c_1(|x|^2+|y|^2)/4}$, and the polynomial factors are absorbed:
\begin{equation*}
 \sup_{0\le\omega\le e^{-1/2}}\sup_{y\in\R^N}e^{\theta'|y|^2/4}
 \int_{\R^N}\bigl|\partial_\omega^j\nabla_y^kK_\omega(y,x)\bigr|
 e^{-\theta''|x|^2/4}\,dx<\infty .
\end{equation*}

(a) By \eqref{eq:gaussid} with $b=1-\omega^2$ and $\gamma=\theta_1$,
\[
 \int_{\R^N}K_\omega(y,x)e^{-\theta_1|x|^2/4}\,dx
 =D^{-N/2}\exp\!\Bigl(-\frac{\theta_1|y|^2}{4D}\Bigr),\qquad
 D=\theta_1+(1-\theta_1)\omega^2\in[\theta_1,1].
\]
As $\theta_1/D\ge\theta_1\ge\theta$, multiplying by $e^{\theta|y|^2/4}$ gives
the bound in (a).

(b) For $0\le s\le1$, (a) with $\theta_1=\theta$ and \eqref{eq:Pbound} give
$$\|\Sop(s)Qf\|_\theta\le\theta^{-N/2}(1+\theta^{-N/2})e^{1/2}e^{-s/2}
\|f\|_\theta.$$ For $s\ge1$ put $\omega=e^{-s/2}\le e^{-1/2}$ and $g=Qf$.
Since $\int g=0$ and $K_0(y,x)=\rho(y)$ does not depend on $x$,
\[
 \Sop(s)g(y)=\int_{\R^N}[K_\omega(y,x)-K_0(y,x)]g(x)\,dx
 =\omega\int_0^1\!\int_{\R^N}
 (\partial_{\omega'}K_{\omega'})(y,x)\big|_{\omega'=r\omega}\,g(x)\,dx\,dr.
\]
By the kernel bound above with $j=1$, $k=0$ and $\theta'=\theta''=\theta$, we have
$$\|\Sop(s)g\|_\theta\le C\omega\|g\|_\theta\le C(1+\theta^{-N/2})
e^{-s/2}\|f\|_\theta.$$ This proves \eqref{eq:Sgapsame}, and \eqref{eq:Sgap}
follows from $\|f\|_\theta\le\|f\|_{\theta_1}$.

(c) $\Lop\partial_{y_1}\rho=-\tfrac12\partial_{y_1}\rho$ and
$\int\partial_{y_1}\rho=0$.

(d) The bounds below justify differentiating under the integral. Let $g=Qf$
and $b=1-e^{-s}$. For $s>0$, $\nabla_yK_\omega=-\frac{y-\omega x}{2b}K_\omega$.
With $\xi=y-\omega x$, $\mu:=1-\theta_1$ and $re^{-\mu r^2}\le(2\mu e)^{-1/2}$,
\[
 |\nabla_yK_\omega(y,x)|
 \le\frac{(2\mu e)^{-1/2}}{\sqrt b}\,\theta_1^{-N/2}
 (4\pi b')^{-N/2}e^{-|\xi|^2/(4b')},\qquad b'=b/\theta_1,
\]
and \eqref{eq:gaussid} with $b'$ and $\gamma=\theta_1$ (so that
$D=\omega^2+b=1$) gives
\[
\int|\nabla_yK_\omega(y,x)|e^{-\theta_1|x|^2/4}\,dx\le
(2\mu e)^{-1/2}\theta_1^{-N/2}b^{-1/2}e^{-\theta_1|y|^2/4}.
\]
As $b\ge(1-e^{-1})s$ for $0<s\le1$, this gives
$\|\nabla\Sop(s)g\|_{\theta_0}\le Cs^{-1/2}\|g\|_{\theta_1}$ there. For $s>1$,
since $\nabla_yK_0=\nabla\rho$ is independent of $x$ and $\int g=0$,
\[
 \nabla\Sop(s)g(y)
 =\omega\int_0^1\!\!\int(\partial_{\omega'}\nabla_yK_{\omega'})\big|_{\omega'=r\omega}\,g(x)\,dx\,dr,
\]
and the kernel bound above with $j=k=1$, $\theta'=\theta_0$,
$\theta''=\theta_1$ yields $\|\nabla\Sop(s)g\|_{\theta_0}\le
Ce^{-s/2}\|g\|_{\theta_1}$. In both ranges
$\|g\|_{\theta_1}\le(1+\theta_1^{-N/2})\|f\|_{\theta_1}$.
\end{proof}

\subsection{Nonlinear decomposition and the ancient orbit}\label{sec:nonlinear}

From now on we fix $0<\theta<1$ and exponents
\begin{equation}\label{eq:weightschoice}
 0<\theta_0<\theta<\theta_1<\min\{p\theta,1\},
\end{equation}
which is possible because $p>1$; constants are allowed to depend on these
weights.
For $f,g\in X_\theta$, the pointwise inequality
$$|\F(r)-\F(s)|\le p(|r|^{p-1}+|s|^{p-1})|r-s|$$ and the margin $\theta_1<p\theta$
(after multiplication by $e^{\theta_1|y|^2/4}$ the factor
$e^{-(p\theta-\theta_1)|y|^2/4}\le1$ remains) give $\F(f)\in X_{\theta_1}$ and
\begin{equation}\label{eq:Fbounds}
\begin{aligned}
 \|\F(f)\|_{\theta_1}&\le\|f\|_\theta^p,\\
 \|\F(f)-\F(g)\|_{\theta_1}
 &\le p\bigl(\|f\|_\theta^{p-1}+\|g\|_\theta^{p-1}\bigr)\|f-g\|_\theta.
\end{aligned}
\end{equation}
By \eqref{eq:Pbound} the same bounds hold for $Q\F$ up to the factor
$C_Q:=1+\theta_1^{-N/2}$. In particular $\F$ maps $C(I;X_\theta)$ into
$C(I;X_{\theta_1})$ for every interval $I$.

We look for an ancient solution of \eqref{eq:veq} on $(-\infty,0]$ which decays
as $\tau\to-\infty$. To this end we split off the neutral mode,
\begin{equation*}
 v(\tau,y)=a(\tau)\rho(y)+z(\tau,y),
 \qquad \int_{\R^N}z(\tau,y)\,dy=0.
\end{equation*}
Projecting the equation onto the neutral and the stable part, we obtain
\begin{equation*}
 a'(\tau)=H(a(\tau),z(\tau)),\qquad
 z_\tau=\Lop z+Q\F(a\rho+z),\qquad
 H(a,z)=\int_{\R^N}\F(a\rho+z)\,dy,
\end{equation*}
and the leading coefficient in the amplitude equation is
\begin{equation*}
 \kappa_N:=\int_{\R^N}\rho^p\,dy
 =(4\pi)^{-Np/2}\Bigl(\frac{4\pi}{p}\Bigr)^{N/2}
 =(4\pi)^{-1}p^{-N/2},
\end{equation*}
using $N(p-1)/2=1$. For $\varepsilon>0$ set
\begin{equation*}
 D_0(\tau):=\varepsilon^{-(p-1)}+(p-1)\kappa_N(-\tau),\qquad
 A_\varepsilon(\tau):=D_0(\tau)^{-1/(p-1)},\qquad \tau\le0 ,
\end{equation*}
This is the ancient solution of $A_\varepsilon'=\kappa_NA_\varepsilon^p$ with
$A_\varepsilon(0)=\varepsilon$. Then $A_\varepsilon$ increases,
$A_\varepsilon\le\varepsilon$, and
\begin{equation}\label{eq:Aasymp}
 A_\varepsilon(\tau)\le[(p-1)\kappa_N(-\tau)]^{-N/2},\qquad
 A_\varepsilon(\tau)\sim C_N(-\tau)^{-N/2},\qquad
 C_N=[\kappa_N(p-1)]^{-N/2}.
\end{equation}
We shall repeatedly use two elementary facts: since $D_0'=-(p-1)\kappa_N$ and
$\log x\le x/\e$,
\begin{equation}\label{eq:Aid}
 \int_\tau^0A_\varepsilon(\sigma)^{p-1}\,d\sigma
 =\frac{\log\bigl(\varepsilon^{p-1}D_0(\tau)\bigr)}{(p-1)\kappa_N},
 \qquad
 \frac{\log\bigl(\varepsilon^{p-1}D_0(\tau)\bigr)}{D_0(\tau)}
 \le\frac{\varepsilon^{p-1}}{\e}\qquad(\tau\le0).
\end{equation}

We solve the equations for $a$ and $z$ simultaneously, by a contraction
argument in a tube around $A_\varepsilon\rho$. Throughout,
$C_{\rm gap}$ and $C_\nabla$ are the constants of \cref{lem:OU}, and
$\|\rho\|_\theta=(4\pi)^{-N/2}\le1$.
Fix $\delta\in(0,1/4)$.  For $\varepsilon>0$ and $M>0$ let
$\mathfrak B_\varepsilon=\mathfrak B_\varepsilon(M)$ be the set of pairs $(a,z)$
with $a\in C((-\infty,0])$ and $z\in C((-\infty,0];X_\theta)$ such that, for every
$\tau\le0$,
\begin{equation}\label{eq:tube}
 (1-\delta)A_\varepsilon(\tau)\le a(\tau)\le(1+\delta)A_\varepsilon(\tau),\qquad
 \int_{\R^N}z(\tau)\,dy=0,\qquad \|z(\tau)\|_\theta\le M A_\varepsilon(\tau)^p .
\end{equation}
For $\eta>0$ equip $\mathfrak B_\varepsilon$ with the metric
$d=d_a+\eta\,d_z$, where
\begin{equation*}
 d_a=\sup_{\tau\le0}\frac{|a_1(\tau)-a_2(\tau)|}{A_\varepsilon(\tau)},\qquad
 d_z=\sup_{\tau\le0}\frac{\|z_1(\tau)-z_2(\tau)\|_\theta}{A_\varepsilon(\tau)^p}.
\end{equation*}
It is easy to see that this metric is complete: a Cauchy sequence converges in
the two weighted
supremum norms, its limit is continuous on compact time intervals, and the
closed tube bounds and mean-zero condition pass to the limit.

\begin{proposition}[Ancient center--stable orbit]\label{prop:LP}
There exists $M_0>0$, depending only on $N,\theta_0,\theta,\theta_1$ and
$\delta$, such that for every $M\ge M_0$ there are $\eta>0$ and
$\varepsilon_0=\varepsilon_0(M)>0$ with the following property.  For every
$0<\varepsilon\le\varepsilon_0$ there is a unique
$(a,z)\in\mathfrak B_\varepsilon(M)$
such that $v=a\rho+z$ solves
\begin{equation}\label{eq:LPsystem}
 a(\tau)^{-(p-1)}=\varepsilon^{-(p-1)}+(p-1)\int_\tau^0B(a(\sigma),z(\sigma))\,d\sigma,
 \qquad
 z(\tau)=\int_{-\infty}^\tau\Sop(\tau-\sigma)Q\F(v(\sigma))\,d\sigma,
\end{equation}
where
\[
B(a,z):=a^{-p}H(a,z)=\int_{\R^N}\F(\rho+z/a)\,dy
\]
and the second integral is defined pointwise in $y$ as in \cref{lem:duhamelinf}.
This solution has the following properties.
\begin{enumerate}[label=\textup{(\roman*)}]
\item $a\in C^1((-\infty,0])$, $a(0)=\varepsilon$, and $a'=a^pB(a,z)>0$; in
particular $v\not\equiv0$.
\item For $\tau\to-\infty$,
\begin{equation}\label{eq:arelerr}
 \frac{a(\tau)}{A_\varepsilon(\tau)}-1
 =O\!\left(\frac{\log(2+|\tau|)}{1+|\tau|}\right),
 \qquad \|z(\tau)\|_\theta\le M A_\varepsilon(\tau)^p .
\end{equation}
\item For all $\tau_0<\tau\le0$ and $y\in\R^N$,
\begin{equation}\label{eq:vmild}
 v(\tau)=\Sop(\tau-\tau_0)v(\tau_0)
 +\int_{\tau_0}^{\tau}\Sop(\tau-\sigma)\F(v(\sigma))\,d\sigma .
\end{equation}
\end{enumerate}
\end{proposition}

\begin{proof}
\textbf{Step 1: constants.}  Put $c_v:=2(4\pi)^{-N/2}+1$,
\[
 C_1:=C_Q\,c_v^p,\qquad
 C_B:=2p\,3^{p-1}\Bigl(\frac{4\pi}{p\theta}\Bigr)^{N/2},\qquad
 C_5:=2p\,c_v^{p-1},\qquad
 C_6:=2C_{\rm gap}C_QC_5 .
\]
Set
\begin{equation*}
 M_0:=2C_{\rm gap}C_1,
\end{equation*}
fix any $M\ge M_0$, put $\eta:=1/(4C_6)$, and finally choose
$\varepsilon_0=\varepsilon_0(M)>0$ so small that for
$0<\varepsilon\le\varepsilon_0$
\begin{equation*}
 \begin{gathered}
 M\varepsilon^{p-1}\le1,\qquad
 2C_BM\varepsilon^{p-1}\le\frac{\kappa_N}2,\qquad
 r_0:=\frac{2C_BM\varepsilon^{p-1}}{\e\,\kappa_N}\le r_\delta,\\
 \varepsilon^{p-1}\Bigl(C_4(1+M)+\frac14\Bigr)\le\min\Bigl\{\frac14,\frac{\eta}{2}\Bigr\},
 \end{gathered}
\end{equation*}
where $C_4:=\beta2^{\beta+2}C_B/(\e\kappa_N)$ and $r_\delta\in(0,1/2]$ is such
that $|(1+r)^{-\beta}-1|\le\delta$ for $|r|\le r_\delta$.  Here
$C_1,C_B,C_4,C_5,C_6$ depend only on $N$, the weights and $\delta$, and in
particular not on $M$ or $\varepsilon$.

\medskip
\textbf{Step 2: bounds in the tube.}  Let $(a,z)\in\mathfrak B_\varepsilon$ and
$v=a\rho+z$.  Since $A_\varepsilon\le\varepsilon$ and $M\varepsilon^{p-1}\le1$,
\begin{equation}\label{eq:vtube}
 \|v(\tau)\|_\theta\le(1+\delta)(4\pi)^{-N/2}A_\varepsilon+MA_\varepsilon^p
 \le c_vA_\varepsilon(\tau),
\end{equation}
and by \eqref{eq:Fbounds}
\begin{equation*}
 \|Q\F(v(\tau))\|_{\theta_1}\le C_1A_\varepsilon(\tau)^p .
\end{equation*}
Let $\zeta=z/a$. Then
\[
\|\zeta\|_\theta\le MA_\varepsilon^p/((1-\delta)A_\varepsilon)
\le2MA_\varepsilon^{p-1}\le2.
\]
Since $\|\rho\|_\theta\le1$, the functions
$\rho$ and $\rho+\zeta$ are bounded by $3e^{-\theta|y|^2/4}$, and the pointwise
inequality for $\F$ gives
\begin{equation}\label{eq:Bclose}
 |B(a,z)-\kappa_N|=\Bigl|\int\bigl(\F(\rho+\zeta)-\F(\rho)\bigr)\Bigr|
 \le C_B\|\zeta\|_\theta\le2C_BMA_\varepsilon^{p-1}.
\end{equation}
In particular $\kappa_N/2\le B\le3\kappa_N/2$ by the smallness conditions of Step~1.  For two
pairs in $\mathfrak B_\varepsilon$, with $\zeta_j=z_j/a_j$,
\begin{equation*}
 \|\zeta_1-\zeta_2\|_\theta
 \le\frac{\|z_1-z_2\|_\theta}{a_1}+\|z_2\|_\theta\frac{|a_1-a_2|}{a_1a_2}
 \le2A_\varepsilon^{p-1}(d_z+Md_a),
\end{equation*}
where we used $(1-\delta)^{-2}<2$, and therefore
\begin{equation*}
 |B(a_1,z_1)-B(a_2,z_2)|\le C_B\|\zeta_1-\zeta_2\|_\theta
 \le2C_B(1+M)A_\varepsilon^{p-1}(d_a+d_z).
\end{equation*}

\medskip
\textbf{Step 3: the Lyapunov--Perron map.}  The identity
\[
\frac d{d\tau}a^{-(p-1)}=-(p-1)a^{-p}a'=-(p-1)B(a,z)
\]
and the terminal value $a(0)=\varepsilon$ suggest
\begin{equation*}
 \Phi_a(a,z)(\tau)=\Bigl[\varepsilon^{-(p-1)}
 +(p-1)\int_\tau^0B(a(\sigma),z(\sigma))\,d\sigma\Bigr]^{-1/(p-1)},
\end{equation*}
which is well defined and continuous because $B>0$ and
$\sigma\mapsto B(a(\sigma),z(\sigma))$ is continuous.  For the mean-zero
component, the condition $z(\tau)\to0$ as $\tau\to-\infty$ suggests
\begin{equation*}
 \Phi_z(a,z)(\tau)=\int_{-\infty}^{\tau}\Sop(\tau-\sigma)Q\F(v(\sigma))\,d\sigma .
\end{equation*}
The map $\sigma\mapsto Q\F(v(\sigma))$ belongs to $C((-\infty,0];X_{\theta_1})$,
has mean zero, and is bounded by $C_1A_\varepsilon^p$ by Step~2.  By
\cref{lem:duhamelinf}, $\Phi_z(a,z)\in C((-\infty,0];X_\theta)$ has mean zero.
Put $\Phi=(\Phi_a,\Phi_z)$.

\medskip
\textbf{Step 4: the center component stays in the tube.}  Write
$\Phi_a(\tau)^{-(p-1)}=D_0(\tau)(1+r(\tau))$ with
$r(\tau)=D_0(\tau)^{-1}(p-1)\int_\tau^0(B-\kappa_N)\,d\sigma$.  By
\eqref{eq:Bclose} and \eqref{eq:Aid},
\begin{equation*}
 |r(\tau)|\le\frac{2C_BM}{\kappa_N}\,
 \frac{\log(\varepsilon^{p-1}D_0(\tau))}{D_0(\tau)}\le r_0\le r_\delta .
\end{equation*}
Since $\Phi_a/A_\varepsilon=(1+r)^{-\beta}$, the choice of $r_\delta$ gives
the first bound in \eqref{eq:tube} for $\Phi_a$.

\medskip
\textbf{Step 5: the stable component stays in the tube.}  By
\cref{lem:duhamelinf}(i), the bound $\|Q\F(v)\|_{\theta_1}\le C_1A_\varepsilon^p$ and the monotonicity of
$A_\varepsilon$,
\begin{equation*}
 \|\Phi_z(\tau)\|_\theta\le C_{\rm gap}C_1\int_{-\infty}^{\tau}
 e^{-(\tau-\sigma)/2}A_\varepsilon(\sigma)^p\,d\sigma
 \le2C_{\rm gap}C_1A_\varepsilon(\tau)^p\leq MA_\varepsilon(\tau)^p .
\end{equation*}
Hence $\Phi(\mathfrak B_\varepsilon)\subset\mathfrak B_\varepsilon$.

\medskip
\textbf{Step 6: contraction in the center direction.}  Let
$(a_j,z_j)\in\mathfrak B_\varepsilon$ and let $D_j$ denote the expression inside
brackets in the definition of $\Phi_a(a_j,z_j)$.  By Step~4,
$D_j\ge D_0(1-r_0)\ge D_0/2$.  By the Lipschitz bound for $B$ in Step~2 and \eqref{eq:Aid},
\[
 |D_1-D_2|\le(p-1)\int_\tau^0|B_1-B_2|\,d\sigma
 \le\frac{2C_B(1+M)}{\kappa_N}(d_a+d_z)\log\bigl(\varepsilon^{p-1}D_0(\tau)\bigr).
\]
By the mean value theorem,
$|D_1^{-\beta}-D_2^{-\beta}|\le\beta(D_0/2)^{-\beta-1}|D_1-D_2|$; dividing by
$A_\varepsilon=D_0^{-\beta}$ and using \eqref{eq:Aid},
\begin{equation*}
 d_a\bigl(\Phi(a_1,z_1),\Phi(a_2,z_2)\bigr)\le C_4(1+M)\varepsilon^{p-1}(d_a+d_z).
\end{equation*}

\medskip
\textbf{Step 7: contraction in the stable direction.}  By \eqref{eq:Fbounds} and
\eqref{eq:vtube},
\[
 \|\F(v_1)-\F(v_2)\|_{\theta_1}
 \le2pc_v^{p-1}A_\varepsilon^{p-1}\bigl(|a_1-a_2|\,\|\rho\|_\theta+\|z_1-z_2\|_\theta\bigr)
 \le C_5A_\varepsilon^p\bigl(d_a+\varepsilon^{p-1}d_z\bigr).
\]
Applying \cref{lem:duhamelinf}(i) to the difference and arguing as in Step~5,
\begin{equation*}
 d_z\bigl(\Phi(a_1,z_1),\Phi(a_2,z_2)\bigr)\le C_6\bigl(d_a+\varepsilon^{p-1}d_z\bigr).
\end{equation*}
Combining the estimates for the two components, with $\eta=1/(4C_6)$, and the
last smallness condition of Step~1 give
\begin{equation*}
 d(\Phi(a_1,z_1),\Phi(a_2,z_2))
 \le\Bigl(C_4(1+M)\varepsilon^{p-1}+\tfrac14\Bigr)d_a
 +\Bigl(C_4(1+M)\varepsilon^{p-1}+\tfrac14\varepsilon^{p-1}\Bigr)d_z
 \le\frac12\bigl(d_a+\eta d_z\bigr).
\end{equation*}
Banach's fixed-point theorem yields a unique fixed point
$(a,z)\in\mathfrak B_\varepsilon$, i.e.\ a solution of \eqref{eq:LPsystem}.

\medskip
\textbf{Step 8: properties.}  (i) The right-hand side of the first equation in
\eqref{eq:LPsystem} is $C^1$ in $\tau$, so $a$ is $C^1$ and
$a'=a^pB(a,z)\ge\kappa_Na^p/2>0$; also $a(0)=\varepsilon$.

(ii) The bound on $z$ is part of the definition of $\mathfrak B_\varepsilon$.
With $r$ as in Step~4 (now for the fixed point), $a/A_\varepsilon-1=(1+r)^{-\beta}-1
=O(|r|)$, and by the first inequality for $|r|$ in Step~4,
$|r(\tau)|\le C\log(\varepsilon^{p-1}D_0(\tau))/D_0(\tau)$ with
$D_0(\tau)\asymp1+|\tau|$.  This is \eqref{eq:arelerr}.

(iii) By construction $\frac d{d\tau}\int v=a'=H(a,z)=\int\F(v)$.  Hence
\[
 a(\tau)\rho=a(\tau_0)\rho+\int_{\tau_0}^\tau P\F(v(\sigma))\,d\sigma,
 \qquad\text{and}\qquad \Sop(\tau-\sigma)P=P.
\]
\Cref{lem:duhamelinf}(iii) gives the corresponding identity for $z$.  Adding
the two identities gives \eqref{eq:vmild}.
\end{proof}

\begin{remark}[Independence of $M$]\label{rem:Mindep}
If $M_0\le M\le M'$ and $0<\varepsilon\le\min\{\varepsilon_0(M),\varepsilon_0(M')\}$,
then $\mathfrak B_\varepsilon(M)\subset\mathfrak B_\varepsilon(M')$ and the metric
does not depend on $M$. The fixed point in $\mathfrak B_\varepsilon(M)$ is
therefore a fixed point in $\mathfrak B_\varepsilon(M')$, and uniqueness there
identifies the two. Thus, for every $\varepsilon$ to which \cref{prop:LP}
applies for some $M\ge M_0$ (we call such $\varepsilon$ \emph{admissible}),
the orbit does not depend on the choice of $M$. We denote it by
$v_\varepsilon=a_\varepsilon\rho+z_\varepsilon$.
\end{remark}

\subsection{Positivity of the ancient orbit}\label{sec:positivity}

The decomposition $v=a\rho+z$ does not show directly that $v$ is positive:
$z$ is controlled only in a weight that is weaker than the decay of $\rho$, so
it might dominate $a\rho$ for large $|y|$. We therefore argue through
approximations on finite time intervals.

\begin{proposition}[Strict positivity of the ancient orbit]\label{prop:positivev}
For every admissible $\varepsilon$,
\[
 v_\varepsilon(\tau,y)>0\qquad(\tau\le0,\ y\in\R^N).
\]
\end{proposition}

\begin{proof}
Fix $M\ge M_0$ and $0<\varepsilon\le\varepsilon_0(M)$, and write $v=v_\varepsilon$.

\medskip
\textbf{Step 1: finite-horizon approximations.} For $n\in\mathbb N$ consider on
$[-n,0]$ the center equation
\begin{equation*}
 a_n(\tau)^{-(p-1)}=\varepsilon^{-(p-1)}
 +(p-1)\int_\tau^0B(a_n(\sigma),z_n(\sigma))\,d\sigma
\end{equation*}
together with the truncated stable equation
\begin{equation*}
 z_n(\tau)=\int_{-n}^{\tau}\Sop(\tau-\sigma)
 Q\F(a_n(\sigma)\rho+z_n(\sigma))\,d\sigma.
\end{equation*}
The proof of \cref{prop:LP} applies verbatim on $[-n,0]$, with the lower limit
$-\infty$ replaced by $-n$ and \cref{lem:duhamelinf}(iv) in place of
\cref{lem:duhamelinf}(i)--(iii). Every stable estimate becomes smaller, while
the center estimates are unchanged. Thus the same constants give a unique pair
$(a_n,z_n)$ satisfying \eqref{eq:tube} on $[-n,0]$, uniformly
in $n$. Put $v_n=a_n\rho+z_n$. Differentiation of the center equation, whose
integrand is continuous, gives $a_n'=H(a_n,z_n)$, and hence
\[
 a_n(\tau)\rho=a_n(-n)\rho
 +\int_{-n}^{\tau}P\F(v_n(\sigma))\,d\sigma.
\]
Adding the truncated stable equation and using $\Sop(s)P=P$, we obtain
\begin{equation}\label{eq:finitehorizonmild}
 v_n(\tau)=\Sop(\tau+n)\bigl[a_n(-n)\rho\bigr]
 +\int_{-n}^{\tau}\Sop(\tau-\sigma)\F(v_n(\sigma))\,d\sigma.
\end{equation}

\medskip
\textbf{Step 2: positivity of $v_n$.} Set $t_n=e^{-n}$ and
$U_n(t,x)=t^{-N/2}v_n(\log t,x/\sqrt t)$ for $t_n\le t\le1$.
By \eqref{eq:conj}, \eqref{eq:finitehorizonmild} is equivalent to
\[
 U_n(t)=e^{(t-t_n)\Delta}U_n(t_n)
 +\int_{t_n}^{t}e^{(t-s)\Delta}\F(U_n(s))\,ds,
 \qquad
 U_n(t_n,x)=t_n^{-N/2}a_n(-n)\rho(x/\sqrt{t_n})>0.
\]
The tube bounds make $U_n$ bounded on $[t_n,1]\times\R^N$. Write
$\F(U_n)=c_nU_n$ with $c_n:=|U_n|^{p-1}$, so that $0\le c_n\le C_n$. Let
$h>0$ satisfy $C_nh<1$. On $[t_n,t_n+h]$, $U_n$ is a fixed point in
$L^\infty$ of the linear map
\[
 \Lambda(W)(t)=e^{(t-t_n)\Delta}U_n(t_n)+\int_{t_n}^te^{(t-s)\Delta}(c_nW)(s)\,ds,
\]
which is a contraction. Since $\Lambda$ preserves nonnegativity, the iterates
starting from $W=0$ are nonnegative, and so is their limit, which is $U_n$.
Restarting at $t_n+h,t_n+2h,\dots$ gives $U_n\ge0$ on $[t_n,1]$. The mild
equation then yields $U_n(t)\ge e^{(t-t_n)\Delta}U_n(t_n)>0$. Returning to
similarity variables, $v_n>0$ on $[-n,0]\times\R^N$.

\medskip
\textbf{Step 3: convergence $v_n\to v$.} The tube gives, uniformly in $n$,
\[
 |v_n(\tau,y)|\le C A_\varepsilon(\tau)e^{-\theta|y|^2/4}
 \qquad(-n\le\tau\le0).
\]
By \cref{app:regularity}, $U_n$ is classical for $t>t_n$. Fix
$k_0\in\mathbb N$ and $\varrho>0$. For $n>k_0+1$ the tube bound gives
\[
 \|U_n\|_{L^\infty([\e^{-k_0-1},1]\times\R^N)}
 \le \e^{N(k_0+1)/2}c_v\varepsilon ,
\]
uniformly in $n$. Restarting the mild identity at $t=\e^{-k_0-1}$, the interior
estimates of \cref{app:regularity} give a bound for $U_n$ in
$C^{2+\gamma,1+\gamma/2}([\e^{-k_0},1]\times \overline{B_\varrho})$ that does
not depend on $n$; the same holds for $v_n$ on $[-k_0,0]\times\overline{B_\varrho}$.
By the Arzel\`a--Ascoli theorem and a diagonal argument in $k_0$ and
$\varrho$, along a subsequence $n_j\to\infty$ we have local uniform
convergence on $(-\infty,0]\times\R^N$ to a classical limit $\widetilde v$.
We set $\widetilde a(\tau):=\int_{\R^N}\widetilde v(\tau,y)\,dy$ and
$\widetilde z:=\widetilde v-\widetilde a\rho$.
Fix $\tau\le0$ and split the stable integral at a fixed $\tau_1<\tau$.  By
\eqref{eq:Sgap} and the tube estimate,
\[
 \left\|\int_{-n_j}^{\tau_1}\Sop(\tau-\sigma)Q\F(v_{n_j}(\sigma))\,d\sigma
 \right\|_\theta
 \le C\int_{-\infty}^{\tau_1}e^{-(\tau-\sigma)/2}
 A_\varepsilon(\sigma)^p\,d\sigma,
\]
and the right-hand side tends to zero as $\tau_1\to-\infty$, uniformly in $j$.
On $[\tau_1,\tau]$ we pass to the limit in the weighted norm:
\begin{equation*}
 \sup_{\sigma\in[\tau_1,\tau]}
 \|\F(v_{n_j}(\sigma))-\F(\widetilde v(\sigma))\|_{\theta_1}
 \longrightarrow0.
\end{equation*}
Indeed, on a fixed ball this follows from local uniform convergence, whereas
outside that ball
\[
 e^{\theta_1|y|^2/4}
 |\F(v_{n_j})-\F(\widetilde v)|
 \le C A_\varepsilon(\sigma)^p
 e^{-(p\theta-\theta_1)|y|^2/4},
\]
which is uniformly small when the radius is large.  Therefore \eqref{eq:Sgap}
and this uniform convergence give dominated convergence in $X_\theta$ on
$[\tau_1,\tau]$, so that
\[
 \int_{-n_j}^{\tau}\Sop(\tau-\sigma)Q\F(v_{n_j}(\sigma))\,d\sigma
 \longrightarrow
 \int_{-\infty}^{\tau}\Sop(\tau-\sigma)Q\F(\widetilde v(\sigma))\,d\sigma
 \qquad\text{in }X_\theta .
\]
By the truncated stable equation, the left-hand side equals
$z_{n_j}(\tau)=v_{n_j}(\tau)-a_{n_j}(\tau)\rho$. The tube bound displayed
above and dominated convergence give, for every fixed $\sigma\le0$,
\[
 a_{n_j}(\sigma)=\int_{\R^N}v_{n_j}(\sigma,y)\,dy
 \longrightarrow \int_{\R^N}\widetilde v(\sigma,y)\,dy=\widetilde a(\sigma),
\]
so that $z_{n_j}(\tau)\to\widetilde v(\tau)-\widetilde a(\tau)\rho
=\widetilde z(\tau)$ pointwise. Comparing the two limits, we obtain
\[
 \widetilde z(\tau)=\int_{-\infty}^{\tau}\Sop(\tau-\sigma)
 Q\F(\widetilde v(\sigma))\,d\sigma.
\]
For the center equation, the lower tube bound for $a_{n_j}$ and the same
Gaussian majorant give an integrable bound for $\F(\rho+z_{n_j}/a_{n_j})$,
uniformly on compact $\sigma$-intervals. Consequently
$B(a_{n_j},z_{n_j})\to B(\widetilde a,\widetilde z)$, and the center equation
passes to the limit. The same majorant gives continuity of $\widetilde a$.
Local uniform convergence and the margin $\theta_1<p\theta$ give continuity of
$\F(\widetilde v)$ into $X_{\theta_1}$, by the argument used for
this uniform convergence, and \cref{lem:duhamelinf} then gives
$\widetilde z\in C((-\infty,0];X_\theta)$. The tube inequalities pass to the
limit. Thus $(\widetilde a,\widetilde z)\in\mathfrak B_\varepsilon(M)$ is a
fixed point of the Lyapunov--Perron map, and uniqueness in \cref{prop:LP} gives
$\widetilde v=v$.  Since every subsequence has the same limit, $v_n\to v$
locally uniformly on $(-\infty,0]\times\R^N$.

\medskip
\textbf{Step 4: strict positivity.} By Steps~2 and~3, $v\ge0$. Moreover
$\int_{\R^N}v(\tau,y)\,dy=a(\tau)>0$, so $v(\tau)\not\equiv0$ for every $\tau$.
For $\tau_0<\tau$, the mild identity \eqref{eq:vmild} becomes
\[
 v(\tau)=\Sop(\tau-\tau_0)v(\tau_0)
 +\int_{\tau_0}^{\tau}\Sop(\tau-\sigma)v(\sigma)^p\,d\sigma
 \ge\Sop(\tau-\tau_0)v(\tau_0)>0,
\]
because the Mehler kernel is strictly positive.
\end{proof}

\subsection{Return to physical variables: proof of \texorpdfstring{\cref{thm:whole}}{Theorem 1.2}}\label{sec:physical}

Fix an admissible $\varepsilon$, let $v=v_\varepsilon=a\rho+z$, and define
\begin{equation*}
 U(t,x)=t^{-N/2}v\!\left(\log t,\frac{x}{\sqrt t}\right),\qquad 0<t\le1 .
\end{equation*}
Write $A(t):=A_\varepsilon(\log t)$.  By \eqref{eq:Aasymp} and since
$L(t)\ge2$ for $t\le T_*$, there are $c,C>0$ with
\begin{equation*}
 cL(t)^{-N/2}\le A(t)\le CL(t)^{-N/2}
 \qquad(0<t\le T_*).
\end{equation*}

\begin{proof}[Proof of \cref{thm:whole}]
\textbf{Step 1: critical estimates.}
The change of variables $x=\sqrt t\,y$ gives $\|U(t)\|_1=\|v(\log t)\|_1$ and
$\int U(t)=a(\log t)$.  By \eqref{eq:tube} and the embedding
$X_\theta\hookrightarrow L^1$,
\[
 (1-\delta)A\le a\le\|v\|_1\le(1+\delta)A+(4\pi/\theta)^{N/2}MA^p\le CA,
\]
and the two-sided bound for $A$ gives the $L^1$ bounds in \eqref{eq:wholebounds}.  Similarly
$\|v\|_\infty\le\|v\|_\theta\le c_vA$ by \eqref{eq:vtube}, which gives
the $L^\infty$ bound in \eqref{eq:wholebounds}.  Then, for $0<t\le T_*$,
\begin{equation*}
 \|\F(U(t))\|_1\le\|U(t)\|_\infty^{p-1}\|U(t)\|_1
 \le\frac{C}{tL(t)^{N/2+1}},
\end{equation*}
which is integrable near $0$ because $N/2+1>1$; on $[T_*,1]$ the same quantity
is bounded. With the $L^1$-isometric dilations $D_th:=t^{-N/2}h(\cdot/\sqrt t)$,
which form a strongly continuous family on $L^1$, we have $U(t)=D_tv(\log t)$
and $\F(U(t))=t^{-1}D_t\F(v(\log t))$. Since $v\in C((-\infty,0];X_\theta)$,
$\F(v)\in C((-\infty,0];X_{\theta_1})$ and
$X_\theta,X_{\theta_1}\hookrightarrow L^1$, both $U$ and $\F(U)$ belong to
$C((0,1];L^1(\R^N))$, and $\F(U)\in L^1((0,1);L^1(\R^N))$.

\medskip
\textbf{Step 2: mild identity from zero.}
Evaluate \eqref{eq:vmild} with $\tau=\log t$, $\tau_0=\log t_0$, $0<t_0<t\le1$,
at $y=x/\sqrt t$ and multiply by $t^{-N/2}$.  By \eqref{eq:conj} with $s=t_0$
and $h=v(\tau_0)$, the first term becomes $e^{(t-t_0)\Delta}U(t_0)$.  For the
integral, put $s=e^\sigma$; since
$\F(U(s))=s^{-1}\cdot s^{-N/2}\F(v(\sigma))(\cdot/\sqrt s)$,
\eqref{eq:conj} gives
\[
 t^{-N/2}\bigl[\Sop(\tau-\sigma)\F(v(\sigma))\bigr](x/\sqrt t)
 =s\,\bigl[e^{(t-s)\Delta}\F(U(s))\bigr](x),
\]
and $d\sigma=ds/s$.  Hence
\begin{equation}\label{eq:Duhamelt0}
 U(t)=e^{(t-t_0)\Delta}U(t_0)+\int_{t_0}^te^{(t-s)\Delta}\F(U(s))\,ds .
\end{equation}
Let $t_0\downarrow0$: the heat semigroup is an $L^1$ contraction, so
$\|e^{(t-t_0)\Delta}U(t_0)\|_1\le\|U(t_0)\|_1\to0$ by \eqref{eq:wholebounds}, and the
integral converges in $L^1$ by Step~1.  Therefore
\begin{equation*}
 U(t)=\int_0^te^{(t-s)\Delta}\F(U(s))\,ds\qquad\text{in }L^1(\R^N),\quad 0<t\le1,
\end{equation*}
and, since $\|U(t)\|_1\to0$, $U\in C([0,1];L^1(\R^N))$ is a mild solution with
datum $0$. Boundedness on $[t_0,1]$ for each $t_0>0$ follows from
$\|v\|_\infty\le c_vA_\varepsilon$, which gives the regularity class asserted in \cref{thm:whole}.

\medskip
\textbf{Step 3: regularity and positivity.} The derivation of
\eqref{eq:Duhamelt0} uses only \eqref{eq:vmild} and \eqref{eq:conj}. By
\cref{app:regularity}, $U$ is a classical solution of \eqref{eq:whole} on
$(0,1]\times\R^N$. \Cref{prop:positivev} gives the positivity of $U$.

\medskip
\textbf{Step 4: the sharp logarithmic endpoint.}
By \eqref{eq:arelerr} and \eqref{eq:Aasymp}, $(-\tau)^{N/2}a(\tau)\to C_N$, while
$(-\tau)^{N/2}\|z(\tau)\|_\theta\le M(-\tau)^{N/2}A_\varepsilon(\tau)^p\to0$.
Hence
\[
 (-\tau)^{N/2}v(\tau,\cdot)\longrightarrow C_N\rho
 \quad\text{in }X_\theta .
\]
Since $X_\theta$ embeds continuously into $L^1$ and $L^\infty$, and
$\|\rho\|_1=1$, $\|\rho\|_\infty=\rho(0)$, this proves \eqref{eq:exactnorms} with
\begin{equation*}
 c_N=\rho(0)\,C_N=(4\pi)^{-N/2}\bigl[(4\pi)^{-1}p^{-N/2}\cdot\tfrac2N\bigr]^{-N/2}
 =\Bigl(\frac N2\Bigr)^{N/2}\Bigl(1+\frac2N\Bigr)^{N^2/4},
\end{equation*}
where we used the formula for $\kappa_N$ in \cref{sec:nonlinear}. For instance $c_1=2^{-1/2}3^{1/4}\approx0.9306$
and $c_2=2$. The first assertion in \eqref{eq:Miy} follows from the $L^\infty$ bound in
\eqref{eq:wholebounds}, and the second from \eqref{eq:exactnorms}.

\medskip
\textbf{Step 5: Gaussian decay.} We show that, with $c:=\theta_0/4$,
\begin{equation}\label{eq:gradUgdecay}
 |U(t,x)|+\sqrt t\,|\nabla U(t,x)|\le Ct^{-N/2}L(t)^{-N/2}e^{-c|x|^2/t}
 \qquad(0<t\le T_*,\ x\in\R^N),
\end{equation}
which contains the Gaussian bound in \eqref{eq:wholebounds}. Since $\rho\le(4\pi)^{-N/2}e^{-\theta|y|^2/4}$,
\eqref{eq:vtube} gives $|v(\tau,y)|\le c_vA_\varepsilon(\tau)e^{-\theta|y|^2/4}$.
For the gradient, by \cref{lem:OU}(d) and \eqref{eq:Fbounds}, the $y$-gradient of
the integrand in
$z(\tau)=\int_{-\infty}^\tau\Sop(\tau-\sigma)Q\F(v(\sigma))\,d\sigma$
is bounded by
$C_\nabla k(\tau-\sigma)c_v^pA_\varepsilon(\sigma)^pe^{-\theta_0|y|^2/4}$,
which is integrable in $\sigma$ uniformly in $y$. Hence we may differentiate
under the integral and
\[
 \|\nabla z(\tau)\|_{\theta_0}\le C_\nabla c_v^p\int_{-\infty}^\tau
 k(\tau-\sigma)A_\varepsilon(\sigma)^p\,d\sigma
 \le C_\nabla c_v^pA_\varepsilon(\tau)^p\int_0^\infty k(s)\,ds
 \le4C_\nabla c_v^pA_\varepsilon(\tau)^p .
\]
Since $|\nabla\rho(y)|=\frac{|y|}2\rho(y)\le Ce^{-\theta_0|y|^2/4}$, we obtain
$|\nabla_yv(\tau,y)|\le CA_\varepsilon(\tau)e^{-\theta_0|y|^2/4}$. Returning to
physical variables and using the two-sided bound for $A$ gives \eqref{eq:gradUgdecay}.
\end{proof}

\begin{remark}[Delayed solutions and the role of $L^\infty_{\rm loc}$]\label{rem:delayed}
For $0<t_1<1$ let $U^{(t_1)}(t)=0$ for $t\le t_1$ and
$U^{(t_1)}(t)=U(t-t_1)$ for $t_1<t\le1$.  By the mild identity from zero in Step~2 of the proof of \cref{thm:whole}, this is a
mild solution with datum $0$ in the sense of \cite[Definition~1.2]{Miyamoto2021},
since $\|U^{(t_1)}(t)\|_\infty<\infty$ for each fixed $t$.  However
$\|U^{(t_1)}(t)\|_\infty\to\infty$ as $t\downarrow t_1$, so $U^{(t_1)}\notin
L^\infty_{\rm loc}((0,1);L^\infty)$.  The requirement $L^\infty_{\rm loc}$ in
the Brezis--Cazenave class excludes such delayed branches; the solutions of
\cref{thm:whole,thm:multicore} belong to that class.
\end{remark}

\begin{proof}[Proof of \cref{cor:scaling}]
For $\lambda\in(0,1]$, critical scaling ($Np=N+2$) shows that
$U_\lambda(t,x)=\lambda^NU(\lambda^2t,\lambda x)$ solves \eqref{eq:whole} on
$(0,1]$, since $\lambda^2t\le1$.  Its $L^1$ trace is still zero, positivity is
preserved, and the endpoint estimates follow from those for $U$, with constants
depending on $\lambda$.  Moreover
\[
 \int_{\R^N}U_\lambda(t,x)\,dx
 =a(\log t+2\log\lambda).
\]
By \cref{prop:LP}(i), $a$ is strictly increasing.  Thus distinct values of
$\lambda$ give distinct solutions.
\end{proof}

\subsection{Rigidity: proof of \texorpdfstring{\cref{thm:rigidity}}{Theorem 1.5}}\label{sec:rigidity}

In this subsection $\theta$ is the weight of \cref{thm:rigidity}, and
$\theta_0,\theta_1$ are chosen as in \eqref{eq:weightschoice}, so that the
results of \cref{sec:nonlinear} apply. The main tool is the spectral gap
\eqref{eq:Sgapsame}, which, unlike \eqref{eq:Sgap}, involves no change of
weight. This is what allows us to treat a solution about which nothing is known
beyond the hypotheses of the theorem.

\begin{proof}[Proof of \cref{thm:rigidity}]
By oddness of the equation we may assume $v\ge0$.

\medskip
\textbf{Step 1: Duhamel formula and mass identity.} For $-\infty<s<\tau<0$, the
function $U(t,x)=t^{-N/2}v(\log t,x/\sqrt t)$ is a bounded classical solution of
$U_t-\Delta U=\F(U)$ on $[\e^s,\e^\tau]\times\R^N$ with Gaussian spatial decay.
The variation-of-constants formula for such solutions \cite{QuittnerSouplet2007},
conjugated through \eqref{eq:conj}, gives
\begin{equation}\label{eq:finiteDuhamel}
 v(\tau)=\Sop(\tau-s)v(s)
 +\int_s^\tau\Sop(\tau-\sigma)\F(v(\sigma))\,d\sigma
\end{equation}
pointwise. By continuity of $v$ into $X_\theta$ and dominated convergence, the
identity extends to $\tau=0$. Since $\Sop$ preserves the integral and
$\F(v)\in C((-\infty,0];X_{\theta_1})\subset C((-\infty,0];L^1)$, integration
over $\R^N$ gives, for $a(\tau):=\int v(\tau,y)\,dy$,
\[
 a(\tau)=a(s)+\int_s^\tau\int_{\R^N}\F(v(\sigma,y))\,dy\,d\sigma,
 \qquad a'(\tau)=\int_{\R^N}\F(v(\tau,y))\,dy .
\]
With $z(\tau):=v(\tau)-a(\tau)\rho$, applying $Q$ to \eqref{eq:finiteDuhamel}
gives
\begin{equation}\label{eq:zDuhamel}
 z(\tau)=\Sop(\tau-s)z(s)
 +\int_s^\tau\Sop(\tau-\sigma)Q\F(v(\sigma))\,d\sigma .
\end{equation}

\medskip
\textbf{Step 2: sign and monotone mass.} Suppose that $v(\tau_2)\equiv0$ for some
$\tau_2$. By Step~1, $U$ is a bounded mild solution on $[\e^{\tau_2},1]$ with
zero datum at $t=\e^{\tau_2}$, and uniqueness of bounded mild solutions
\cite{Weissler1980,QuittnerSouplet2007} gives $v\equiv0$ on $[\tau_2,0]$. If $v(\tau)\not\equiv0$
for some $\tau<\tau_2$, then \eqref{eq:finiteDuhamel} and $\F(v)\ge0$ give
$v(\tau_2)\ge\Sop(\tau_2-\tau)v(\tau)>0$, a contradiction. Hence $v\equiv0$,
contrary to the hypothesis. Therefore $v(\tau)\not\equiv0$ for every $\tau$, and
\eqref{eq:finiteDuhamel} gives $v(\tau)\ge\Sop(\tau-s)v(s)>0$. In particular
$v>0$ on $(-\infty,0]\times\R^N$ and $a'=\int v^p>0$.

\medskip
\textbf{Step 3: center dominance.} Since $X_\theta\hookrightarrow L^1$,
$\|z(s)\|_\theta\le C\|v(s)\|_\theta\to0$ as $s\to-\infty$. Moreover
$\|v^p\|_\theta\le\|v\|_\theta^p$ and $Q$ is bounded on $X_\theta$. By
\eqref{eq:Sgapsame}, each term of \eqref{eq:zDuhamel} is bounded pointwise by
its $X_\theta$-norm times $e^{-\theta|y|^2/4}$, so
\[
 e^{\theta|y|^2/4}|z(\tau,y)|\le C_\theta e^{-(\tau-s)/2}\|z(s)\|_\theta
 +C\int_s^\tau e^{-(\tau-\sigma)/2}\|v(\sigma)\|_\theta^p\,d\sigma .
\]
Letting $s\to-\infty$ gives
\begin{equation*}
 \|z(\tau)\|_\theta\le C\int_{-\infty}^{\tau}
 e^{-(\tau-\sigma)/2}\|v(\sigma)\|_\theta^p\,d\sigma,
\end{equation*}
and, since $\sup_{\sigma\le0}\|v(\sigma)\|_\theta<\infty$, the same limit in
\eqref{eq:zDuhamel} gives the pointwise identity
\begin{equation*}
 z(\tau)=\int_{-\infty}^\tau\Sop(\tau-\sigma)Q\F(v(\sigma))\,d\sigma .
\end{equation*}
Let $h(\tau):=\sup_{\sigma\le\tau}\|v(\sigma)\|_\theta$. By the
estimate for $\|z(\tau)\|_\theta$ just obtained, $\|z(\sigma)\|_\theta\le Ch(\tau)^p$ for $\sigma\le\tau$,
and since $a$ is increasing,
$\|v(\sigma)\|_\theta\le a(\sigma)\|\rho\|_\theta+\|z(\sigma)\|_\theta
\le Ca(\tau)+Ch(\tau)^p$. Hence $h(\tau)\le Ca(\tau)+Ch(\tau)^p$, and since
$h(\tau)\to0$ the last term can be absorbed. For all sufficiently negative
$\tau$, with $C$ depending only on $N$ and $\theta$,
\begin{equation}\label{eq:autocenter}
 \sup_{\sigma\le\tau}\|v(\sigma)\|_\theta\le C a(\tau),
 \qquad
 \|z(\tau)\|_\theta\le C a(\tau)^p .
\end{equation}

\medskip
\textbf{Step 4: asymptotics.} Since $a>0$, $a'=a^pB(a,z)$ with $B$ as in
\cref{prop:LP}. As $\|z/a\|_\theta=O(a^{p-1})$, the estimate \eqref{eq:Bclose}
gives
\begin{equation*}
 a'(\tau)=\kappa_Na(\tau)^p+O\bigl(a(\tau)^{2p-1}\bigr)
 =a(\tau)^p\left[\kappa_N+O(a(\tau)^{p-1})\right].
\end{equation*}
In particular, for sufficiently negative $\tau$,
$\frac{\kappa_N}{2}a^p\le a'\le2\kappa_Na^p$, so
$a(\tau)^{p-1}\asymp(-\tau)^{-1}$.  Equivalently,
$\frac d{d\tau}a^{-(p-1)}=-(p-1)\kappa_N+O(a^{p-1})$, and integration from a
fixed negative time to $\tau$ gives
$a(\tau)^{-(p-1)}=(p-1)\kappa_N(-\tau)+O(\log(-\tau))$.
This proves the asymptotics of $a$ in (i), and together with \eqref{eq:autocenter}
it yields the profile convergence in (i). This completes part (i).

\medskip
\textbf{Step 5: identification with the Lyapunov--Perron orbit.} Let $\tau_1<0$
be so negative that \eqref{eq:autocenter}, the bound $\|z/a\|_\theta\le2$ and
the two-sided inequality of Step~4 hold on $(-\infty,\tau_1]$, and let
$\tau_*\le\tau_1$. Put $\varepsilon:=a(\tau_*)$, translate $\tau_*$ to
zero, and denote the translated quantities again by $a,z$. With
$D_0(s)=\varepsilon^{-(p-1)}+(p-1)\kappa_N(-s)$ and
$A_\varepsilon=D_0^{-1/(p-1)}$, the amplitude equation of Step~4 gives
\[
 a(s)^{-(p-1)}=\varepsilon^{-(p-1)}
 +(p-1)\int_s^0\left[\kappa_N+O(a(r)^{p-1})\right]dr .
\]
The lower bound $a'\ge\frac{\kappa_N}2a^p$ implies
$a(r)^{p-1}\le C/(\varepsilon^{-(p-1)}+c(-r))$, whence
$\int_s^0a(r)^{p-1}dr\le C\log(1+\varepsilon^{p-1}|s|)$. Dividing the resulting
denominator error by $D_0(s)$ and using
$\sup_{x\ge0}\log(1+x)/(1+x)<\infty$ gives
\begin{equation}\label{eq:entrycenter}
 \sup_{s\le0}\left|\frac{a(s)}{A_\varepsilon(s)}-1\right|
 \le C\varepsilon^{p-1}.
\end{equation}
Then \eqref{eq:autocenter} and \eqref{eq:entrycenter} give
\begin{equation}\label{eq:entrystable}
 \|z(s)\|_\theta\le C A_\varepsilon(s)^p\qquad(s\le0).
\end{equation}
The constants in \eqref{eq:entrycenter}--\eqref{eq:entrystable} depend on
$\tau_1$ but not on $\tau_*\le\tau_1$. Fix $M\ge\max\{M_0,C\}$ with $C$ from
\eqref{eq:entrystable}. There is $\bar\tau\le\tau_1$ such that, for
$\tau_*\le\bar\tau$, $\varepsilon=a(\tau_*)\le
\varepsilon_0(M)$ and $C\varepsilon^{p-1}\le\delta$ in \eqref{eq:entrycenter}.
For such $\tau_*$ the translated pair $(a,z)$ lies in $\mathfrak B_\varepsilon(M)$:
$a$ is continuous, $z\in C((-\infty,0];X_\theta)$ has mean zero, and the tube
bounds are \eqref{eq:entrycenter}--\eqref{eq:entrystable}. It is a fixed point
of the Lyapunov--Perron map: the center equation follows by integrating
$\frac d{ds}a^{-(p-1)}=-(p-1)B(a,z)$ from $s$ to $0$ with $a(0)=\varepsilon$,
and the stable equation is the identity for $z$ obtained in Step~3. Uniqueness in \cref{prop:LP} and
\cref{rem:Mindep} give $v(\cdot+\tau_*)=v_\varepsilon$ on $(-\infty,0]$. This
proves part (ii).

\medskip
\textbf{Step 6: uniqueness up to translation.} Let $\bar\varepsilon$ be
admissible and $\bar v:=v_{\bar\varepsilon}$, with mass $\bar a$. By
\cref{prop:LP}, \cref{prop:positivev} and Step~3 of the proof of
\cref{thm:whole}, $\bar v$ satisfies the hypotheses of \cref{thm:rigidity}. Let
$\bar\tau$ and $\bar s$ be the thresholds of Step~5 for $v$ and $\bar v$. Since
$a$ and $\bar a$ are continuous, strictly increasing and tend to $0$ at
$-\infty$, every $\varepsilon\le\min\{a(\bar\tau),\bar a(\bar s)\}$ is attained
as $\varepsilon=a(\tau_*)=\bar a(s_*)$ with $\tau_*\le\bar\tau$ and
$s_*\le\bar s$. By Step~5 and \cref{rem:Mindep},
$v(\cdot+\tau_*)=v_\varepsilon=\bar v(\cdot+s_*)$ on $(-\infty,0]$. With
$s:=s_*-\tau_*$ this gives $v(\tau)=\bar v(\tau+s)$ for $\tau\le\tau_*$, and
uniqueness of bounded mild solutions \cite{Weissler1980,QuittnerSouplet2007},
applied in physical variables from the time $\e^{\tau_*}$, extends the identity
to $\tau\le\min\{0,-s\}$. This proves part (iii).
\end{proof}

\section{Localization over arbitrary bounded backgrounds}\label{sec:localization}

In this section we prove \cref{thm:multicore} by placing finitely many
truncated copies of the whole-space solution inside $\Omega$.

\subsection{Cutoff cores and their defects}\label{sec:cutoffcores}

Throughout this section the constants may depend on $N$, $\Omega$, the
whole-space solution $U$, and the configuration of cores and cutoff functions
(in particular on the distances between the cores and their distance to
$\partial\Omega$). We write $C_{K,m}$ when we want to stress the dependence on
the bound $K$ of the background and on the number $m$ of cores. These
constants are not uniform as cores merge or approach the boundary; for a single
core in a compact subset of $\Omega$ uniformity is proved at the end of the
section.

Fix distinct points $x_1,\dots,x_m\in\Omega$.  Choose radii $r_j>0$ so small
that the closed balls $\overline{B_{2r_j}(x_j)}$ are pairwise disjoint and
contained in $\Omega$.  Choose $\chi_j\in C_c^\infty(B_{2r_j}(x_j))$ with
$0\le\chi_j\le1$ and $\chi_j=1$ on $B_{r_j}(x_j)$.  Let $U$ be the positive
whole-space core of \cref{thm:whole} and set
\begin{equation*}
 W_j(t,x):=\chi_j(x)U(t,x-x_j).
\end{equation*}
Its defect is
\begin{equation*}
 E_j:=(W_j)_t-\Delta W_j-\F(W_j)
 =(\chi_j-\chi_j^p)\F(U)-2\nabla\chi_j\cdot\nabla U-(\Delta\chi_j)U,
\end{equation*}
where $U$ is evaluated at $(t,x-x_j)$. The defect is supported where
$r_j\le|x-x_j|\le2r_j$. The Gaussian bound \eqref{eq:gradUgdecay} for $U$ and
$\nabla U$, together with $|\F(U)|=|U|^p$ and $p-1=2/N$, gives, after
decreasing $c$,
\begin{equation}\label{eq:cutoffbounds}
 |E_j(t,x)|\le Ct^{-(N+2)/2}e^{-c/t},\qquad
 |W_j(t,x)|^{p-1}\le \frac{C}{tL(t)}
 e^{-c|x-x_j|^2/t}\qquad(0<t\le T_*,\ x\in\Omega).
\end{equation}

For a sign vector $\boldsymbol\sigma\in\{-1,+1\}^m$, put
\begin{equation*}
 W_{\boldsymbol\sigma}:=\sum_{j=1}^m\sigma_jW_j,
 \qquad E_{\boldsymbol\sigma}:=\sum_{j=1}^m\sigma_jE_j.
\end{equation*}
Because the cutoff supports are disjoint,
$\F(W_{\boldsymbol\sigma})=\sum_{j}\sigma_j\F(W_j)$, and therefore
\begin{equation*}
 (W_{\boldsymbol\sigma})_t-\Delta W_{\boldsymbol\sigma}
 =\F(W_{\boldsymbol\sigma})+E_{\boldsymbol\sigma}.
\end{equation*}
Since $W_{\boldsymbol\sigma}$ is smooth on $(0,T_*]\times\overline\Omega$ and
vanishes near $\partial\Omega$, for $0<t_0<t\le T_*$ the linear
variation-of-constants formula gives
\[
 W_{\boldsymbol\sigma}(t)=S_\Omega(t-t_0)W_{\boldsymbol\sigma}(t_0)
 +\int_{t_0}^tS_\Omega(t-s)
 [\F(W_{\boldsymbol\sigma}(s))+E_{\boldsymbol\sigma}(s)]\,ds .
\]
By \eqref{eq:wholebounds}, $\|W_{\boldsymbol\sigma}(t_0)\|_1\le C_mL(t_0)^{-N/2}\to0$
as $t_0\downarrow0$. Moreover, \eqref{eq:wholebounds} and \eqref{eq:cutoffbounds}
give $\|\F(W_{\boldsymbol\sigma}(s))\|_1\le
\|W_{\boldsymbol\sigma}(s)\|_\infty^{p-1}\|W_{\boldsymbol\sigma}(s)\|_1
\le C_ms^{-1}L(s)^{-N/2-1}$ and, after decreasing $c$,
$\|E_{\boldsymbol\sigma}(s)\|_1\le C_me^{-c/s}$; both bounds are integrable at
$s=0$. Letting $t_0\downarrow0$ and using the $L^1$ contraction of $S_\Omega$,
we obtain
\begin{equation*}
 W_{\boldsymbol\sigma}(t)
 =\int_0^tS_\Omega(t-s)
 [\F(W_{\boldsymbol\sigma}(s))+E_{\boldsymbol\sigma}(s)]\,ds
\end{equation*}
in $L^1(\Omega)$. The identity also holds pointwise: for $0\le s\le t/2$ the
kernel of $S_\Omega(t-s)$ is bounded by $Ct^{-N/2}$, so, once $t_0\le t/2$, the
two terms that vanish in $L^1$ as $t_0\downarrow0$, namely
$S_\Omega(t-t_0)W_{\boldsymbol\sigma}(t_0)$ and the integral over $[0,t_0]$,
also vanish uniformly in $x$.

\subsection{The bounded background and the localized Duhamel gain}

Fix $u_0\in L^\infty(\Omega)$ and put $K=\max\{2\|u_0\|_\infty,1\}$.
The standard bounded-data theory
\cite{Weissler1980,BrezisCazenave1996,QuittnerSouplet2007} gives, for
sufficiently small $T_0>0$, a unique bounded mild solution
\begin{equation*}
 V(t)=S_\Omega(t)u_0+\int_0^tS_\Omega(t-s)\F(V(s))\,ds,
\end{equation*}
with
\begin{equation*}
 \|V\|_{L^\infty((0,T_0)\times\Omega)}\le K.
\end{equation*}
It is classical for positive time, satisfies the Dirichlet condition and
converges to $u_0$ in $L^1$. If $u_0\ge0$, comparison gives $V\ge0$.

The logarithmic gain comes from the following elementary estimate.

\begin{lemma}[Localized critical Duhamel estimate]\label{lem:localizedD}
Let $c_0>0$, $\mu>0$, and $x_0\in\Omega$.  There is $C=C(N,c_0,\mu)$ such that,
for $0<t\le T_*$,
\begin{equation*}
 \left\|\int_0^tS_\Omega(t-s)
 \left[\frac{e^{-c_0|\cdot-x_0|^2/s}}{sL(s)^\mu}\right]ds\right\|_\infty
 \le \frac{C}{L(t)^\mu}.
\end{equation*}
If $N=1$, then also
\begin{equation*}
 \left\|\int_0^tS_\Omega(t-s)
 \left[s^{-1/2}L(s)^{-1/2}e^{-c_0|\cdot-x_0|^2/s}\right]ds\right\|_\infty
 \le Ct^{1/2}L(t)^{-1/2}.
\end{equation*}
\end{lemma}

\begin{proof}
The Dirichlet heat kernel is dominated by the Euclidean kernel.  Writing
$e^{-c_0|x|^2/s}$ as a multiple of the heat kernel at time $s/(4c_0)$ gives,
for $0<s<t$,
\[
 \left\|e^{(t-s)\Delta}
 e^{-c_0|\cdot-x_0|^2/s}\right\|_\infty
 =\Bigl(\frac{s}{4c_0(t-s)+s}\Bigr)^{N/2}
 \le\max\{1,(4c_0)^{-N/2}\}\left(\frac{s}{t}\right)^{N/2}.
\]
Hence
\[
 \left\|\int_0^tS_\Omega(t-s)
 \frac{e^{-c_0|\cdot-x_0|^2/s}}{sL(s)^\mu}\,ds\right\|_\infty
 \le Ct^{-N/2}\int_0^ts^{N/2-1}L(s)^{-\mu}\,ds
 \le CL(t)^{-\mu}.
\]
When $N=1$, the same computation gives the second estimate.
\end{proof}

We shall use the elementary inequalities
\begin{equation*}
 |\F(A+B)-\F(A)-\F(B)|\le C_p|A|\,|B|^{p-1},
 \qquad 1<p\le2,
\end{equation*}
which follows by integrating $A[\F'(B+rA)-\F'(rA)]$ over $0\le r\le1$ and using
the $(p-1)$-H\"older continuity of $\F'(s)=p|s|^{p-1}$, and, for $p=3$,
\begin{equation*}
 (A+B)^3-A^3-B^3=3AB(A+B).
\end{equation*}

Put $Z_{\boldsymbol\sigma}:=V+W_{\boldsymbol\sigma}$ and define its defect
\begin{align*}
 R_{\boldsymbol\sigma}
 &:=(Z_{\boldsymbol\sigma})_t-\Delta Z_{\boldsymbol\sigma}
   -\F(Z_{\boldsymbol\sigma})
 =E_{\boldsymbol\sigma}+\F(V)+\F(W_{\boldsymbol\sigma})
   -\F(V+W_{\boldsymbol\sigma}).
\end{align*}
If $N\ge2$, then $1<p\le2$, so the first elementary inequality, the bound $|V|\le K$,
and \eqref{eq:cutoffbounds} give
\begin{equation*}
 |R_{\boldsymbol\sigma}(t,x)|
 \le |E_{\boldsymbol\sigma}(t,x)|
 +\frac{CK}{tL(t)}\sum_{j=1}^m e^{-c|x-x_j|^2/t}.
\end{equation*}
For $N=1$, the identity for $p=3$ gives
\begin{equation*}
 |R_{\boldsymbol\sigma}|
 \le |E_{\boldsymbol\sigma}|+C(K|W_{\boldsymbol\sigma}|^2+K^2|W_{\boldsymbol\sigma}|).
\end{equation*}
Consequently \cref{lem:localizedD} and \eqref{eq:cutoffbounds} imply
\begin{equation*}
 \sup_{0<t\le T}L(t)
 \left\|\int_0^tS_\Omega(t-s)R_{\boldsymbol\sigma}(s)\,ds\right\|_\infty
 \le C_{K,m}
\end{equation*}
for all sufficiently small $T$.

\subsection{The correction and the bounded-domain theorem}

Let
\begin{equation*}
 \mathcal X_T:=\left\{w\in C((0,T];L^\infty(\Omega)):\
 \|w\|_{\mathcal X_T}:=\sup_{0<t\le T}L(t)\|w(t)\|_\infty<\infty\right\}.
\end{equation*}
Equipped with the displayed norm, $\mathcal X_T$ is a Banach space: a Cauchy
sequence converges uniformly in $L^\infty(\Omega)$ on each interval
$[t_0,T]$, $t_0>0$, and the weighted bound passes to the limit.

\begin{proposition}[Correction over a bounded background]\label{prop:logcorrection}
For every fixed finite configuration $(\mathbf x,\boldsymbol\sigma)$ there is
$T\in(0,\min\{T_0,T_*\}]$ and a unique fixed point in a suitable ball of
$\mathcal X_T$ of
\begin{equation}\label{eq:Psilog}
 \Psi_{\boldsymbol\sigma}(w)(t)
 :=\int_0^tS_\Omega(t-s)
 \left[\F(Z_{\boldsymbol\sigma}+w)-\F(Z_{\boldsymbol\sigma})
       -R_{\boldsymbol\sigma}\right](s)\,ds.
\end{equation}
It satisfies
\begin{equation}\label{eq:wlog}
 \|w(t)\|_\infty\le\frac{C_{K,m}}{L(t)}
\end{equation}
and
\begin{equation}\label{eq:wL1sharp}
 \|w(t)\|_1\le C_{K,m}\left(\frac{t^{N/2}}{L(t)}
 +\frac{t}{L(t)^p}+e^{-c/t}\right)
 =o(L(t)^{-N/2}).
\end{equation}
\end{proposition}

\begin{proof}
For real $A,B_1,B_2$,
\[
 |\F(A+B_1)-\F(A+B_2)|
 \le C_p\bigl(|A|^{p-1}+|B_1|^{p-1}+|B_2|^{p-1}\bigr)|B_1-B_2|.
\]
Since the cutoff supports are disjoint, \eqref{eq:cutoffbounds} gives
\[
 |W_{\boldsymbol\sigma}|^{p-1}
 \le\frac{C}{tL(t)}\sum_{j=1}^m e^{-c|x-x_j|^2/t},
\]
and hence
$|Z_{\boldsymbol\sigma}|^{p-1}\le C(K^{p-1}+|W_{\boldsymbol\sigma}|^{p-1})$.
If $\|w_i\|_{\mathcal X_T}\le r_*$ and
$d=\|w_1-w_2\|_{\mathcal X_T}$, the bounded-background term contributes
$CK^{p-1}Td$.  By \cref{lem:localizedD} with $\mu=2$, the singular coefficient
contributes at most $C_mL(T)^{-1}d$, and the pure-correction term contributes
$Cr_*^{p-1}TL(T)^{1-p}d$.  Thus
\begin{equation*}
 \|\Psi_{\boldsymbol\sigma}(w_1)-\Psi_{\boldsymbol\sigma}(w_2)\|_{\mathcal X_T}
 \le C\left(L(T)^{-1}+K^{p-1}T+r_*^{p-1}TL(T)^{1-p}\right)d.
\end{equation*}
We still have to check that $\Psi_{\boldsymbol\sigma}(w)$ belongs to
$\mathcal X_T$.
Write $G_w=\F(Z_{\boldsymbol\sigma}+w)-\F(Z_{\boldsymbol\sigma})
-R_{\boldsymbol\sigma}$. The preceding pointwise bounds, integrated in
space, give $G_w\in L^1((0,T);L^1(\Omega))$; it is bounded on each
strip $[t_0,T]\times\Omega$, $t_0>0$. Near zero the core terms are bounded in
$L^1$ by multiples of $s^{N/2-1}/L(s)$ and $s^{N/2-1}/L(s)^2$. In dimension
one the additional term has size $O(L(s)^{-1/2})$, which is also integrable.

Fix $t_1>0$ and choose $0<t_2<t_1$. For $t$ near $t_1$, the integral over
$[0,t_2]$ is continuous into $C_0(\overline\Omega)$ by the positive time
gap, heat-kernel smoothing and dominated convergence in $L^1$.
For the integral over $[t_2,t]$, remove a terminal interval of length
$h$. Its $L^\infty$ norm is at most
$h\sup_{t_2\le s\le T}\|G_w(s)\|_\infty$, whereas the remaining
integral is continuous by the same smoothing argument. Letting
$h\downarrow0$ proves continuity into $L^\infty$ and the zero boundary
values for positive time. (We argue in this way because $S_\Omega$ is not
strongly continuous on $L^\infty$.)

Together with the bound for the Duhamel integral of $R_{\boldsymbol\sigma}$ in the previous subsection, this shows that, with
$r_*=2C_{K,m}$ and $T$ sufficiently small, $\Psi_{\boldsymbol\sigma}$ is a
contraction of the closed ball of radius $r_*$ into itself.  This proves
\eqref{eq:wlog}.

For the $L^1$ estimate, use the $L^1$ contraction of $S_\Omega$ in
\eqref{eq:Psilog}. By the Lipschitz bound above,
\[
 \|\F(Z_{\boldsymbol\sigma}+w)-\F(Z_{\boldsymbol\sigma})\|_1
 \le CK^{p-1}\|w\|_1
 +C\bigl\|\bigl(|W_{\boldsymbol\sigma}|^{p-1}+|w|^{p-1}\bigr)|w|\bigr\|_1 .
\]
The Gaussian bound gives, in every dimension,
\[
 \||W_{\boldsymbol\sigma}(s)|^{p-1}\|_1
 \le C_m\frac{s^{N/2-1}}{L(s)},
\]
so, by \eqref{eq:wlog}, the singular coefficient contributes at most
$C_mt^{N/2}/L(t)^2$ after integration in time, and the pure-correction term at
most $C_{K,m}t/L(t)^p$. For the defect, if $N\ge2$, then
$\|R_{\boldsymbol\sigma}(s)\|_1\le Ce^{-c/s}+C_{K,m}s^{N/2-1}/L(s)$, whose time
integral is at most $C_{K,m}t^{N/2}/L(t)+Ce^{-c/t}$. If $N=1$, then
the one-dimensional bound for $R_{\boldsymbol\sigma}$ above, $\|W_{\boldsymbol\sigma}(s)\|_2^2\le Cs^{-1/2}L(s)^{-1}$ and
$\|W_{\boldsymbol\sigma}(s)\|_1\le CL(s)^{-1/2}$ give
\[
 \|R_{\boldsymbol\sigma}(s)\|_1\le Ce^{-c/s}
 +C\bigl(Ks^{-1/2}L(s)^{-1}+K^2L(s)^{-1/2}\bigr).
\]
The time integrals of the last two terms are $O(t^{1/2}/L(t))$ and
$O(tL(t)^{-1/2})$, and the latter is bounded by $Ct^{1/2}/L(t)$ for small $t$
because $t^{1/2}L(t)^{1/2}\to0$. Combining these estimates yields
\begin{align*}
 \|w(t)\|_1
 &\le CK^{p-1}\int_0^t\|w(s)\|_1\,ds
 +C_{K,m}\frac{t^{N/2}}{L(t)}
 +C_m\frac{t^{N/2}}{L(t)^2}
 +C_{K,m}\frac{t}{L(t)^p}+Ce^{-c/t}.
\end{align*}
Gronwall's inequality gives \eqref{eq:wL1sharp} after changing the constants.

We note that \eqref{eq:wlog} alone gives only
$\|w(t)\|_1\le|\Omega|C_{K,m}L(t)^{-1}$. This is not $o(L(t)^{-N/2})$ when
$N\ge2$, and not even $O(L(t)^{-N/2})$ when $N\ge3$. The separate $L^1$
estimate \eqref{eq:wL1sharp} shows that the correction does not alter the
mass scale $L(t)^{-N/2}$ of the cores, and it gives the $L^1$ bound in
\cref{thm:multicore}(ii).
\end{proof}

\begin{proof}[Proof of \cref{thm:multicore}]
Let $V$ be the bounded solution constructed above, let
$W_{\boldsymbol\sigma}$ be the signed sum of cores of \cref{sec:cutoffcores}, and let $w$ be the
correction from \cref{prop:logcorrection}.  Set
\begin{equation*}
 u_{\mathbf x,\boldsymbol\sigma}:=V+W_{\boldsymbol\sigma}+w.
\end{equation*}
Adding the mild identities for $V$ and $W_{\boldsymbol\sigma}$ to the fixed
point equation \eqref{eq:Psilog} gives
\[
 u_{\mathbf x,\boldsymbol\sigma}(t)
 =S_\Omega(t)u_0+\int_0^tS_\Omega(t-s)
 \F(u_{\mathbf x,\boldsymbol\sigma}(s))\,ds.
\]
We check the integrability required in \eqref{eq:mildOmega}.
Since $V$ is bounded, \eqref{eq:wlog} gives $\|w(t)\|_\infty\le C/L(t)$, and the
Gaussian estimate for each core yields
\[
 \|W_{\boldsymbol\sigma}(t)\|_p^p
 \le C_m t^{-Np/2}L(t)^{-Np/2}t^{N/2}
 =\frac{C_m}{tL(t)^{N/2+1}},
\]
where we used $N(p-1)/2=1$ and $Np/2=N/2+1$.  Hence
\[
 \|\F(u_{\mathbf x,\boldsymbol\sigma}(t))\|_1
 =\|u_{\mathbf x,\boldsymbol\sigma}(t)\|_p^p
 \le C\left(1+\frac{1}{tL(t)^{N/2+1}}+\frac{1}{L(t)^p}\right),
\]
which is integrable at $t=0$ because $N/2+1>1$.  Thus
$\F(u_{\mathbf x,\boldsymbol\sigma})\in L^1((0,T);L^1(\Omega))$.
The construction gives boundedness on every strip $[t_0,T]\times\Omega$,
$t_0>0$, and continuous zero boundary values. Restarting the mild equation at
time $t_0/2$, the regularity argument of \cref{app:regularity} shows that
$u_{\mathbf x,\boldsymbol\sigma}$ is classical on $[t_0,T]\times\Omega$ and
continuous up to the lateral boundary, with zero boundary values.  Since
$W_{\boldsymbol\sigma}(t)\to0$ in $L^1$ and \eqref{eq:wL1sharp} tends to zero,
$u_{\mathbf x,\boldsymbol\sigma}(t)\to u_0$ in $L^1$.

The $L^\infty$ estimate in (ii) follows from the Gaussian bound
for the finitely many cores and \eqref{eq:wlog}; the $L^1$ estimate in (ii) follows
from \eqref{eq:wholebounds} and \eqref{eq:wL1sharp}.  At $x_j$, the cutoff balls are
disjoint, so $W_j(t,x_j)=U(t,0)$ and $W_k(t,x_j)=0$ for $k\ne j$.  Therefore
\[
 u_{\mathbf x,\boldsymbol\sigma}(t,x_j)-V(t,x_j)
 =\sigma_jU(t,0)+w(t,x_j),
\]
and \eqref{eq:mastercenter} follows from the last limit in \eqref{eq:exactnorms} because
$w(t,x_j)=O(L(t)^{-1})$ is negligible on the core scale.  If $u_0\ge0$, then
$V\ge0$ is bounded; hence any core with $\sigma_j=-1$ forces the solution to be
negative at $x_j$ for all sufficiently small $t$.

Finally let $\mathcal K\Subset\Omega$ and consider only one core.  Choose a
single cutoff profile with radius depending only on
$\operatorname{dist}(\mathcal K,\partial\Omega)$.  All constants above are then
uniform for $x\in\mathcal K$, so the existence time can be chosen independently
of $x$ and of the sign.  If $x\ne y$, then at the point $x$ the core centered
at $y$ is $O(t^{-N/2}L(t)^{-N/2}e^{-c|x-y|^2/t})$ and its correction is
$O(L(t)^{-1})$, whereas the branch centered at $x$ has the nonzero limit
\eqref{eq:mastercenter}. For the same point and opposite signs,
\eqref{eq:mastercenter} gives values of opposite sign relative to $V$. Thus all
single-core branches are pairwise distinct.
\end{proof}

\section{Further remarks and comparison with earlier work}\label{sec:literature}

\Cref{tab:comparison} summarizes the comparison carried out in this section.
Apart from the other doubly critical endpoint $p=N/(N-2)$, the earlier
constructions of nonuniqueness for regular or zero data live in Lebesgue
spaces strictly below the critical exponent $q_c=N(p-1)/2$. At $p=p_F$ this
range is empty, because $q_c=1$. Constructions based on heteroclinic
connections or singular stationary states
\cite{FilaMatanoYanagida2017,HisaMiyamoto2026} start from singular initial
data; they are discussed in \cref{sec:selfsimilar,sec:singular}.

\begin{table}[htbp]
\centering
\small
\renewcommand{\arraystretch}{1.15}
\begin{tabular}{@{}>{\raggedright\arraybackslash}p{4.3cm}
 >{\raggedright\arraybackslash}p{3.4cm}
 >{\raggedright\arraybackslash}p{2.7cm}
 >{\raggedright\arraybackslash}p{4.4cm}@{}}
\toprule
Work & Exponent & Trace space & Initial data\\
\midrule
Haraux--Weissler, Weissler
\cite{HarauxWeissler1982,Weissler1986a,Weissler1986b}
& $p_F<p<p_S$ & $L^q$, $1\le q<q_c$ & $u_0=0$ on $\R^N$ (positive self-similar profile)\\
Glogi\'c--Hofmanov\'a--Lange--Luongo \cite{GlogicHofmanovaLangeLuongo2025}
& $N\ge3$, $p_F<p<p_{JL}$ & $L^q$, $1\le q<q_c$ & specific data in $L^q(\R^N)$\\
Cazenave--Dickstein--Naumkin--Weissler
\cite{CazenaveDicksteinNaumkinWeissler2020,CazenaveDicksteinNaumkinWeissler2019}
& a range of $p>1$ & $L^r$, $1\le r<q_c$ & $\mu|x|^{-2/(p-1)}$; perturbations on bounded domains\\
Ni--Sacks, Terraneo, Matos--Terraneo
\cite{NiSacks1985,Terraneo2002,MatosTerraneo2003}
& $N\ge3$, $p=N/(N-2)$ & $L^{q_c}=L^{N/(N-2)}$ & every datum in $L^{N/(N-2)}(\R^N)$ \cite{MatosTerraneo2003}\\
\midrule
This paper
& $p=p_F$ & $L^{q_c}=L^1$ & $u_0=0$ on $\R^N$; every $u_0\in L^\infty(\Omega)$\\
\bottomrule
\end{tabular}
\caption{Nonuniqueness mechanisms for $u_t-\Delta u=|u|^{p-1}u$ and the
space in which the initial trace is taken. Here $q_c=N(p-1)/2$, $p_S$ is the
Sobolev exponent and $p_{JL}$ the Joseph--Lundgren exponent.}
\label{tab:comparison}
\end{table}

\subsection{The two doubly critical endpoints}

Miyamoto's uniqueness theorem covers the classes with $\ell>N/2$ recalled in
\cref{sec:intro}, and the case $\ell=N/2$ was left open in
\cite[Remark~1.6]{Miyamoto2021}. The positive solution of \cref{thm:whole} and
the zero solution have the same initial trace in $L^1$ and both belong to the
endpoint class, while the asymptotics \eqref{eq:exactnorms} show that the
positive solution belongs to no class with $\ell>N/2$. We note that the
existence result in $X_{N/2}$ and the nonexistence examples below this exponent
describe sharpness within the logarithmic Orlicz scale; they do not give a
necessary condition for an individual datum.

There is another doubly critical endpoint, $p=q=N/(N-2)$ with $N\ge3$, for which
nonuniqueness was proved by Ni and Sacks \cite{NiSacks1985} and by Terraneo
\cite{Terraneo2002}. Matos and Terraneo \cite{MatosTerraneo2003} showed that
every datum in $L^{N/(N-2)}(\R^N)$ gives rise to infinitely many solutions
which are continuous with values in that space. \Cref{thm:multicore} concerns
the other endpoint, $q_c=1$ and $p=1+2/N$.

\subsection{Self-similar and connecting-orbit mechanisms}\label{sec:selfsimilar}

Positive self-similar solutions with zero initial trace were first
constructed by Haraux and Weissler \cite{HarauxWeissler1982}, and Weissler
\cite{Weissler1986a,Weissler1986b} found rapidly decaying profiles in the
energy-subcritical range
\[
 p_F<p<p_S,\qquad
 p_S:=\begin{cases}
 (N+2)/(N-2),&N\ge3,\\
 +\infty,&N=1,2.
 \end{cases}
\]
Related nonuniqueness results on bounded domains are due to Baras
\cite{Baras1983}, and the profile theory was developed further by Escobedo and
Kavian \cite{EscobedoKavian1987}, Yanagida \cite{Yanagida1996} and Dohmen and
Hirose \cite{DohmenHirose1998}. These solutions have zero trace in $L^q$ for
$1\le q<N(p-1)/2$, and the recent construction of Glogi\'c, Hofmanov\'a, Lange
and Luongo \cite{GlogicHofmanovaLangeLuongo2025} (nonuniqueness in $L^q(\R^N)$
for $N\ge3$, $p_F<p<p_{JL}$ and $1\le q<q_c$) also lives in this subcritical
range. At $p=p_F$ the range is empty. Indeed, if $u(t,x)=t^{-N/2}f(x/\sqrt t)$
with $f\in L^1$, then $\|u(t)\|_1=\|f\|_1$ for all $t$, so a nontrivial
self-similar solution cannot have zero trace in $L^1$; this is why the factor
$L(t)^{-N/2}$ is needed.

Fila and Yanagida \cite{FilaYanagida2011} constructed homoclinic and
heteroclinic orbits in the Sobolev critical and supercritical regimes, and Fila,
Matano and Yanagida \cite{FilaMatanoYanagida2017} used heteroclinic connections
and ancient solutions of a rescaled equation to obtain nonuniqueness without
self-similarity for singular initial data. Their point of view is close to
ours, but the initial traces and the range of parameters are different.

\subsection{Singular data and localization}\label{sec:singular}

Naito \cite{Naito2004,Naito2006} developed a multiplicity theory for positive
self-similar solutions with homogeneous singular data, and Souplet and Weissler
\cite{SoupletWeissler2003} proved multiplicity for data above the singular
steady state. Cazenave, Dickstein, Naumkin and Weissler
\cite{CazenaveDicksteinNaumkinWeissler2020,CazenaveDicksteinNaumkinWeissler2019}
constructed sign-changing self-similar solutions of $u_t-\Delta u=|u|^{p-1}u$,
for a range of $p>1$, with datum $\mu|x|^{-2/(p-1)}$, together with
perturbations on bounded domains. Their perturbation results
\cite{CazenaveDicksteinNaumkinWeissler2019} allow data equal to
$\mu|x-x_0|^{-2/(p-1)}$ near $x_0$, including $\mu=0$, so their Remark~6.3
already contains a construction for bounded data that vanish near a point.
However, the strong continuity in $L^r$ on bounded domains obtained there holds
for $1\le r<N(p-1)/2$, a range which is empty at $p=p_F$; moreover, at this
exponent a nonzero integrable self-similar profile has constant $L^1$ norm. In
our construction the core has zero strong $L^1$ trace, and the background near
each concentration point is an arbitrary bounded function. At the critical
Sobolev exponent, Cort\'azar, del Pino and Musso
\cite{CortazarDelPinoMusso2020} constructed solutions on bounded domains that
bubble in infinite time at several points, which must satisfy a positivity
condition on a matrix built from the Dirichlet Green's function; in
\cref{thm:multicore} the cores have Gaussian tails, and their positions and
signs are arbitrary.

Nonuniqueness can also come from singular stationary states; see the recent
work of Hisa and Miyamoto \cite{HisaMiyamoto2026} and, for exponential
nonlinearities in two dimensions, of Fujishima, Ioku, Ruf and Terraneo
\cite{FujishimaIokuRufTerraneo2026}. Our solutions are of a different nature:
they are bounded and classical at each positive time, and their singular
behavior consists in the concentration of the solution, with blow-up of the
supremum norm, as $t\downarrow0$.

\subsection{Scope of the rigidity result}

\Cref{thm:rigidity} only concerns ancient solutions of fixed sign that tend
to zero in $X_\theta$. Within this class it identifies a single nonstationary
orbit of the autonomous rescaled equation, unique up to time translation,
rather than a stationary self-similar profile. Two ingredients force this: the
spectral gap in a fixed weight and the positivity of the derivative of the
mass. We do not know what happens outside the Gaussian weighted class.

\subsection{Positivity and mass}\label{sec:posmass}

For each $t_0>0$ the usual maximum principle holds on $[t_0,T]\times\Omega$.
This is not in conflict with the loss of positivity at $t=0$: convergence in
$L^1$ gives no control of a negative part that concentrates in $L^\infty$. For
the positive whole-space solution the mass satisfies
\[
 \int_{\R^N}U(t,x)\,dx=a(\log t),\qquad
 \frac{d}{dt}\int_{\R^N}U(t,x)\,dx
 =\int_{\R^N}U(t,x)^p\,dx
 \sim\frac{\kappa_N a(\log t)^p}{t}.
\]
The derivative is integrable near $t=0$, so the mass balance holds with zero
initial mass, even though the solution is nontrivial at every positive time.

\subsection{The absorption problem and necessary conditions on traces}\label{sec:absorption}

The amplitude law of \cref{prop:LP} has a well-known counterpart. For the
absorption equation $u_t=\Delta u-u^p$ with $p=p_F$, Galaktionov, Kurdyumov and
Samarskii \cite{GalaktionovKurdyumovSamarskii1985} showed that there are
nonnegative solutions with $L^1$ data whose large-time behavior is described by
\[
 c_N\bigl[(T+t)\log(T+t)\bigr]^{-N/2}\exp\!\Bigl(-\frac{|x|^2}{4(T+t)}\Bigr),
\]
with exactly the constant $c_N$ of \cref{thm:whole}; see also Herraiz
\cite{Herraiz1999} for a detailed asymptotic analysis of absorption problems.
In similarity variables the absorption amplitude obeys, to leading order,
$a'=-\kappa_Na^p$ as $\tau\to+\infty$, whereas here $a'=\kappa_Na^p$ as
$\tau\to-\infty$, which explains why the same constant appears.

The logarithmic rate is also sharp from the point of view of necessary
conditions. For nonnegative solutions at $p=p_F$, Hisa and Ishige
\cite[Theorem~1.1 and Lemma~3.3]{HisaIshige2018} proved that, for all
sufficiently small radii $\varrho$, the mass in balls of radius $\varrho$ at
times $0<s<\varrho^2$ is bounded by $C[\log(\e+\sqrt T/\varrho)]^{-N/2}$, where
$T$ is the existence time. For the
solution of \cref{thm:whole}, the change of variables $x=\sqrt t\,y$ and the
tube bounds give
\[
 \int_{B_{\sqrt{2t}}(0)}U(t,x)\,dx
 =a(\log t)\int_{B_{\sqrt2}(0)}\rho\,dy+\int_{B_{\sqrt2}(0)}z(\log t,y)\,dy
 \ge cL(t)^{-N/2}
\]
for small $t$. Thus, for $\varrho=\sqrt{2t}$, $U$ attains the Hisa--Ishige
bound up to constants. In particular this necessary condition is sharp, and
the rate $L(t)^{-N/2}$ in \eqref{eq:wholebounds} is the largest possible for
the mass of nonnegative solutions in parabolic balls.

\appendix
\section{Auxiliary estimates}
\subsection{A Gaussian identity}\label{app:gaussian}

Completing the square gives the following identity.
Let $\omega\ge0$, $b>0$ and $\gamma>0$.  Then for $y\in\R^N$
\begin{equation}\label{eq:gaussid}
 (4\pi b)^{-N/2}\int_{\R^N}
 \exp\!\Bigl(-\frac{|y-\omega x|^2}{4b}-\frac{\gamma|x|^2}{4}\Bigr)dx
 =D^{-N/2}\exp\!\Bigl(-\frac{\gamma|y|^2}{4D}\Bigr),\qquad D=\omega^2+\gamma b .
\end{equation}
In particular, for $b=1-\omega^2$, $\omega\in[0,1)$, one has
$D=\gamma+(1-\gamma)\omega^2$.

\subsection{Duhamel integrals from \texorpdfstring{$-\infty$}{minus infinity}}\label{app:functional}

\begin{lemma}\label{lem:duhamelinf}
Let $0<\theta<\theta_1<1$ and let $g\in C((-\infty,0];X_{\theta_1})$ satisfy
$\int g(\sigma)=0$ for every $\sigma$ and $G:=\sup_\sigma\|g(\sigma)\|_{\theta_1}<\infty$.
For $\tau\le0$ and $y\in\R^N$ define
\[
 Z(\tau,y):=\int_{-\infty}^\tau\bigl[\Sop(\tau-\sigma)g(\sigma)\bigr](y)\,d\sigma .
\]
Then:
\begin{enumerate}[label=\textup{(\roman*)}]
\item $Z(\tau)\in X_\theta$, $\int Z(\tau)=0$ and
$\|Z(\tau)\|_\theta\le C_{\rm gap}\int_{-\infty}^\tau e^{-(\tau-\sigma)/2}\|g(\sigma)\|_{\theta_1}\,d\sigma$;
\item $Z\in C((-\infty,0];X_\theta)$;
\item $Z(\tau)=\Sop(\tau-\tau_0)Z(\tau_0)+\int_{\tau_0}^\tau\Sop(\tau-\sigma)g(\sigma)\,d\sigma$
pointwise, for $\tau_0<\tau\le0$;
\item the same conclusions hold for $g\in C([-n,0];X_{\theta_1})$, with $-\infty$
replaced by $-n$ and $\tau,\tau_0\in[-n,0]$.
\end{enumerate}
\end{lemma}

\begin{proof}
For fixed $y$, the map $\sigma\mapsto[\Sop(\tau-\sigma)g(\sigma)](y)$ is continuous
on $(-\infty,\tau)$ by dominated convergence, since the kernel is continuous in
$s=\tau-\sigma>0$ and $g$ is continuous into $X_{\theta_1}$.  Since $Qg=g$,
\eqref{eq:Sgap} gives
\begin{equation*}
 \bigl|[\Sop(\tau-\sigma)g(\sigma)](y)\bigr|\le C_{\rm gap}e^{-(\tau-\sigma)/2}
 \|g(\sigma)\|_{\theta_1}e^{-\theta|y|^2/4},
\end{equation*}
so the integral converges absolutely and the bound in (i) holds.  Continuity of
$Z(\tau)$ in $y$ follows by dominated convergence, and $\int Z(\tau)=0$ by
Fubini's theorem, since $\Sop$ preserves the integral.

For (ii), let $h>0$ with $\tau+h\le0$.  Substituting $\sigma\mapsto\sigma+h$ in
$Z(\tau+h)$ gives
\[
 Z(\tau+h)-Z(\tau)=\int_{-\infty}^\tau\Sop(\tau-\sigma)\bigl[g(\sigma+h)-g(\sigma)\bigr]d\sigma,
\]
hence, for any $\Lambda>0$,
\[
 \|Z(\tau+h)-Z(\tau)\|_\theta
 \le 2C_{\rm gap}\sup_{\tau-\Lambda\le\sigma\le\tau}\|g(\sigma+h)-g(\sigma)\|_{\theta_1}
 +4C_{\rm gap}Ge^{-\Lambda/2}.
\]
Choosing $\Lambda$ large and then $h$ small, by uniform continuity of $g$ on
$[\tau-\Lambda,0]$, proves right continuity.  Applying the same identity at the
point $\tau-h$ proves left continuity.

For (iii), by the semigroup property and Fubini's theorem (justified by the pointwise
bound in the proof of (i)),
\[
\Sop(\tau-\tau_0)Z(\tau_0)=\int_{-\infty}^{\tau_0}\Sop(\tau-\sigma)g(\sigma)\,d\sigma.
\]

For (iv), the arguments for (i) and (iii) are unchanged. In (ii), the
substitution produces, besides the integral over $[-n,\tau]$, an additional
integral over an interval of length $h$, which is bounded in $X_\theta$ by
$C_{\rm gap}Gh$.
\end{proof}

\subsection{Regularity away from the initial time}\label{app:regularity}

We use standard heat-equation regularity; see \cite[Chapter~IV]{LSU1968} and
\cite{QuittnerSouplet2007}. A mild solution of $u_t-\Delta u=\F(u)$ that is
bounded on every strip $[t_0,T]\times\R^N$ (or $[t_0,T]\times\Omega$),
$t_0>0$, is classical away from the initial time. Indeed, restarting the mild
identity at any $t_0>0$, local $W^{2,1}_r$ estimates with $r>N+2$ and the
parabolic Sobolev embedding give H\"older continuity; the local Lipschitz
continuity of $\F$ then gives a H\"older continuous source, and Schauder
estimates yield $C^{2+\gamma,1+\gamma/2}$ regularity locally, for some
$0<\gamma<1$. For $U$ the restarted identity is \eqref{eq:Duhamelt0}, and for
the finite-horizon solutions it is the physical form of
\eqref{eq:finitehorizonmild}, which applies for $t>t_n$.

On a bounded $C^{2,\alpha}$ domain, $S_\Omega(t)u_0$ is continuous on
$\overline\Omega$, vanishes on the boundary for $t>0$, and converges to $u_0$
in $L^1$ as $t\downarrow0$. If a source is integrable in $L^1$ near zero and
bounded on every positive-time strip, its Duhamel integral is continuous up to
the boundary with zero boundary values; this follows by splitting off the
early-time integral, where the kernel has a positive time gap. The same
bootstrap applies up to the lateral boundary, with boundary Schauder estimates
away from the initial time. No compatibility between $u_0$ and the boundary
values at $t=0$ is required, since the initial trace is taken in $L^1$.

\section*{Acknowledgments}
The author is grateful to Professor Manuel Del Pino (UK) for valuable comments and discussions, which helped to improve the paper.

\section*{Declaration on the use of generative AI}

The author used Claude (Anthropic) for language editing, improving the presentation of the manuscript, and assisting in the verification of mathematical arguments. All AI-assisted material was carefully reviewed and independently verified by the author. The author takes full responsibility for the content, correctness, and integrity of the manuscript.

\section*{Declaration of competing interest}
The author declare that there is no conflict of interest.

\section*{Data Availability Statements}
The manuscript has no associated data.

\end{document}